\documentclass[12pt]{article}
\usepackage[T2A]{fontenc}
\usepackage[cp1251]{inputenc}
\usepackage[russian,english]{babel}
\usepackage{amsthm}
\usepackage[intlimits]{amsmath}
\usepackage{amsfonts,amssymb}
\usepackage{graphicx}  
\usepackage{icomma}
\usepackage{bbm}
\usepackage{tikz}
\usetikzlibrary{shapes,snakes}
\usepackage{color}
\usepackage{algorithmic}
\usepackage{verbatim}
\usepackage{xurl}
\usepackage{subfig}
\usepackage{tempora}
\usepackage[upint]{newtxmath}
\usepackage[explicit]{titlesec}
\usepackage{totcount}
\usepackage{caption}
\usepackage{indentfirst}
\usepackage{layout}
\usepackage[a4paper, textwidth=165mm, textheight=235mm, left=22.5mm, top=26mm, footskip=15mm]{geometry}
\usepackage{appendix}

\newlength{\interwordspace}
\newlength{\thinwordspace}
\titleformat{\section}
{\large\centering\normalfont}{\thesection.}{\interwordspace}{\MakeUppercase{#1}}

\titleformat{\subsection}
{\normalsize\centering\normalfont\itshape}{\thesubsection.}{\interwordspace}{#1}

\titleformat{\subsubsection}[runin]
{\normalsize\normalfont\bfseries}{\indent\thesubsubsection.}{\thinwordspace}{#1.}

\newtotcounter{citnum}
\def\oldbibitem{} \let\oldbibitem=\bibitem
\def\bibitem{\stepcounter{citnum}\oldbibitem}
\regtotcounter{table}
\regtotcounter{figure}

\makeatletter
\newcommand\Appendix{
\setcounter{section}{0}%
\setcounter{subsection}{0}%
\gdef\@secapp{\appendixname}%
\gdef\thesection{\@asbuk\c@section}}
\makeatother

\makeatletter
\renewenvironment{abstract}{%
	\small%
	\list{}{\listparindent 0em%
		\itemindent \listparindent
		\rightmargin 10mm \leftmargin 10mm
		\parsep \z@ \@plus\p@}%
	\item\relax}%
{\endlist}
\makeatother

\newcommand{\keywords}[1]{%
\medskip		
{{\noindent\normalsize\bfseries Key words:}~{\small #1}.}%
}

\newcommand{\DOI}[1]{%
\medskip	
	{{\noindent\normalsize\bfseries DOI:}~{\small #1}}%
}

\numberwithin{equation}{section}

\theoremstyle{plain}

\theoremstyle{definition}



\makeatletter
\renewcommand*{\@biblabel}[1]{\hfill#1.}
\makeatother

\let\alpha\alphaup
\let\beta\betaup
\let\gamma\gammaup
\let\delta\deltaup
\let\epsilon\epsilonup
\let\zeta\zetaup
\let\eta\etaup
\let\theta\thetaup
\let\iota\iotaup
\let\kappa\kappaup
\let\lambda\lambdaup
\let\mu\muup
\let\nu\nuup
\let\xi\xiup
\let\pi\piup
\let\rho\rhoup
\let\sigma\sigmaup
\let\tau\tauup
\let\upsilon\upsilonup
\let\phi\phiup
\let\chi\chiup
\let\psi\psiup
\let\omega\omegaup
\let\varepsilon\varepsilonup
\let\vartheta\varthetaup
\let\varpi\varpiup
\let\varrho\varrhoup
\let\varsigma\varsigmaup
\let\varphi\varphiup
\let\varkappa\varkappaup

\usepackage{cite,hyperref}

\newcommand{\AB}{\mathbf{A}}
\newcommand{\FB}{\mathbf{F}}
\newcommand{\IB}{\mathbf{I}}
\newcommand{\cc}{\mathbf{c}}
\newcommand{\ff}{\mathbf{f}}
\newcommand{\yy}{\mathbf{y}}
\newcommand{\Rr}{\mathbb{R}}

\begin{document}
	
\vspace*{0mm}

\begin{center} 
\large\bf
\uppercase{Efficient and stable time integration of Cahn-Hilliard equations: 
explicit, implicit and explicit iterative schemes}
\end{center}

\begin{center}
\textbf{\large\copyright\ 2024\quad  
  M.A.~Botchev,\textsuperscript{1,*}, I.A.~Fahurdinov, \textsuperscript{1,2,**},
  E.B.~Savenkov\textsuperscript{1,***}}\\

{\it
\textsuperscript{1} 
Keldysh Institute of Applied Mathematics of Russian Academy of Sciences,\\
Miusskaya Sq.~4, Moscow, 125047, Russia\\
\textsuperscript{2} National Research Nuclear University 
``Moscow Engineering Physics Institute'' (MEPhI)\\ 
Kashirskoe Hw.~31, Moscow, 115409, Russia\\
*e--mail: botchev@kiam.ru\\
** mv1451003@gmail.com \\
*** savenkov@keldysh.ru
} \\

Submitted\quad --.--.2024

Revised version\quad --.--.2024

Accepted for publication\quad --.--.2024
\end{center} 

\begin{abstract}

To solve the Cahn-Hilliard equation numerically, 
a new time integration algorithm is proposed, which is 
based on a combination of the Eyre splitting and 
the local iteration modified (LIM) scheme.  The latter is employed
to tackle the implicit system arising each time integration step.
The proposed method is gradient-stable and allows to use large time steps,
whereas, regarding its computational structure, it is an explicit 
time integration scheme.
Numerical tests are presented to demonstrate abilities of the new method 
and to compare it with other time integration methods for Cahn-Hilliard equation.
Citations:~\total{citnum}. Figures:~\total{figure}.

\keywords{Cahn-Hilliard equation, gradient-stable schemes, Eyre splitting,
local iteration modified scheme, LIM}

\DOI{}
\end{abstract}


\section{Introduction}
The Cahn-Hilliard equation was proposed in~1958 (see original paper~\cite{cahn_1958})
to describe a phase separation in two-component alloys. 
Currently the equation is widely applied in various research fields
and is one of the fundamental models in the so called gradient or
weakly non-local thermomechanics of continuous media~\cite{gurtin_1996}. 
As a part of more complex models, the equation is employed in
multiphase hydrodynamics, material science, solidification problems,
phase transition theory and in many other fields~\cite{provatas_2010, steinbach_2023}.
From a theoretical point of view, the Cahn-Hilliard equation
forms the basis of phenomenological spinodal decomposition theory~\cite{skripov_1979}
and has a broad range of applications in theoretical physics, 
see~\cite{hohenberg_1977,penrose_1990,bray_2002} and references therein.
Mathematics behind the Cahn-Hilliard equation is also widely discussed
in the literature, see, e.g.,~\cite{miranville_2019}.

Up to recently, the Cahn-Hilliard equation had been primarily seen
as a theoretical tool in mathematical modeling of various processes.
However, last decennia it has been actively employed also as a simulation tool,
either independently or a part of more complex models.
This has necessitated development of efficient numerical algorithms 
to solve the equation.

Difficulties arising in design of efficient numerical algorithms for
the Cahn-Hilliard equation are mainly caused by the following
two reasons.
First, the Cahn-Hilliard equation is a nonlinear PDE containing fourth
order space derivatives.  Due to nonlinearities, solutions of the Cahn-Hilliard 
equation evolve on a broad range of space and time scales.
In particular, a typical solution is almost constant within the certain
space regions of homogeneity which correspond to ``pure'' phases of the system.
These homogeneity regions have a typical size~$d$ and are separated from each 
other by thin layers, of width $\epsilon \ll d$, the so-called 
``diffuse interfaces'', where the solution varies 
from its minimal to maximal values but remains smooth.
At initial stages of spinodal decomposition,
usually emerging due to a perturbation in the initial random distribution,
the diffuse interface width~$\epsilon$ is comparable to homogeneity region
size~$d$ ($\epsilon \sim d$) and
typical evolution times are of order~$\epsilon^2$.
However, in the course of further evolution, when the diffuse interface width
becomes considerably smaller than the homogeneity region size,
$\epsilon \ll d$,
typical evolution times are of order~$1/\epsilon$, 
see~\cite{pego_1989,bates_1993}.

As a consequence of these effects, explicit time integration schemes,
being applied to the Cahn-Hilliard equation, 
are stable for time step sizes~$\tau \sim h^4$, with $h$ being the
space grid size.
Although such a restriction is physically motivated, in some cases
it may turn out unacceptably strict for a computationally efficient solution process.
At the same time, in fully implicit schemes, 
uniqueness of the solution on the next time level is not guaranteed
for sufficiently large time step sizes~\cite{demello_2004,vollmayr_2003}.
Therefore, design of numerical schemes which would be
(a)~uniquely solvable for arbitrary time steps,
(b)~stable, and (c)~conservative,
is not an easy task. 
Baseline algorithms meeting these requirements appeared relatively 
recently~\cite{eyre_1997,eyre_1998}.
They are based on a special splitting of system energy called
convex splitting, see Section~\ref{sec:algo}. 
We note that, in the context of numerical solution of 
the Cahn-Hilliard equation, ``stability'' usually means the so-called 
``gradient stability'' (also known as ``energy stability''),
which guarantees that a discretized analogue of the system free energy
does not grow with time, see 
Sections~\ref{sec:model},~\ref{sec:algo} below.
This stability condition can be considered as strong in the
sense that the expression for the free energy includes first derivatives
of solution. This is necessary to make sure that the numerical solution
of the Cahn-Hilliard equation is thermodynamically consistent.

Currently, important research directions in efficient
time integration of the Cahn-Hilliard equation include
construction of adaptive numerical algorithms (needed due to
the broad range of characteristic evolution times,
as discussed above), 
development of energy-stable time integration schemes
(producing thermodynamically consistent numerical solutions),
and design of schemes which would be, on one hand,
sufficiently asymptotically stable and, on the other hand,
uniquely solvable and computationally efficient.

A complete review of current research on design and analysis of 
time integration schemes for the Cahn-Hilliard equation is far beyond
the scope of this paper.
An overview of main research directions in this field is given in
a survey~\cite{tierra_2013}.
As examples of other relevant literature, 
we note~\cite{felgueroso_2008,li_2017,zhang_2012,minkoff_2006,luo_2016,gonzalez_2014,kassam_2005} 
(adaptive and high-order methods),
\cite{he_2007,song_2015,li_2022} (schemes, allowing ``large'' time step sizes),
and 
\cite{eyre_1997,eyre_1998,chen_2019a,chen_2019b,zhang_2020,zhou_2023,lee_2023,boyer_2011,brachet_2020}
(energy-stable schemes).

To discretize the Cahn-Hilliard equation in space, various methods can
be employed, such as classic finite difference methods, finite volume methods,
spectral methods, standard and isogeometric finite element methods,
see~\cite{elliott_1989, barret_1995, chen_1996, furihata_2001, feng_2004, wells_2006, wise_2009, du_2011, xia_2007, brenner_2020, gomez_2009, zhang_2019, kastner_2016, goudenege_2012}.

The aim of this paper is to study numerically a new class of time integration
schemes for the Cahn-Hilliard equation.  The new algorithms are based
on two main ideas: 
(a)~employment of implicit-explicit time integration 
schemes based on the convex splitting of the system 
energy~\cite{eyre_1997,eyre_1998}, ensuring the scheme is energy-stable, and
(b)~application of the local iteration modified (LIM) time integration scheme,
which allows to combine the computational simplicity of explicit
schemes with stability properties of implicit schemes.
The Eyre splitting~\cite{eyre_1997,eyre_1998} can currently be seen
as a basis for constructing time integration schemes for the
Cahn-Hilliard equations, which should be both energy stable
and uniquely solvable (meaning that the solution on the next time level
in the time integration scheme exists and is unique). 

Local iteration schemes form a family of explicit time integration 
methods, where stability is achieved by employing Chebyshev polynomials.
Chebyshev polynomial iterations are used in time integration since at least 
the 1950s~\cite{CzaoDin1957, CzaoDin1960}.
It should be emphasized that Chebyshev iteration in the local iteration schemes
are used not to approximately solve linear systems arising in implicit
time integration; if this would be the case then, as shown
in~1952 by Gel'fand and Lokutsievskii 
(see~\cite{GelLok1962} and \cite[Ch.10, \S~4.12]{Babenko2002}),
no essential gain in computational work with respect to explicit schemes
can be obtained.
The key characteristic feature of the local iteration schemes is that
Chebyshev iterations are constructed to ensure stability and accuracy
of time integration, rather than to achieve a fast convergence
to a solution of an implicit scheme. 
The local iteration schemes are proposed 
in~\cite{LokLok1984, Zhukov1984, LokLok1986, Zhukov1986}
and further developed in~\cite{Zhukov1993, ShvedZhuk1998, zhukov_2011},
see also references in~\cite{zhukov_2011}. 
Specially tuned Chebyshev iterations allow the local iteration schemes
to work with time step sizes essentially larger than in explicit Euler
scheme, while keeping the computational simplicity of explicit schemes.

The first local iteration scheme proposed in~\cite{Zhukov1984, LokLok1984, LokLok1986}
is not monotone (its numerical solution is not guaranteed to be nonnegative)
and is not asymptotically stable for large times $t\rightarrow\infty$.
Therefore here the local iteration modified (LIM) scheme, proposed in~\cite{Zhukov1986}, 
is employed, which does have there properties.
A detailed description and comparison of the local iteration schemes 
can be found in~\cite{zhukov_2011}. 
Up to date, the LIM scheme has been successfully applied to solve 
various rather challenging real-life problems, 
see, e.g.,~\cite{zhukov_2014,zhukov_2019,zhukov_2021}.
In~\cite{zhukov_2012}, the LIM scheme is applied to solve 
the Cahn-Hilliard equation within a numerical implementation
of a full mathematical model to describe metal crystallization
processes.
There, at each time step, the Cahn-Hilliard equation is linearized 
and the linearized problem is handled by the LIM method.
No comparisons with other relevant time integration schemes
are carried out in~\cite{zhukov_2012}, and properties such as
gradient stability are not investigated either.

In this paper, with a space one-dimensional Cahn-Hilliard equation
taken as example, we propose new time integration schemes 
which are based the LIM scheme and implicit-explicit Eyre splittings.
Results of numerical experiments, where the proposed methods are tested 
and compared to other time integration schemes, are presented.

In Section~\ref{sec:model} a short description of the model is given.
Currently used time integration schemes are briefly discussed in
Section~\ref{sec:algo} and then, in the same section,
new schemes are presented.
Numerical experiments are described in Section~\ref{sec:num}.
Finally, Section~\ref{sec:conclusions} summarizes results
obtained in this paper.

Work of E.B.~Savenkov is supported by 
by the Russian Science Foundation grant no.~2-11-00203.

\section{Mathematical model}\label{sec:model}
We consider an one-dimensional (1D) Cahn-Hilliard equation
\begin{equation}\label{eq:ch}
  \frac{\partial c}{\partial t}
  = \frac{\partial}{\partial x}\left( M\frac{\partial \mu}{\partial x}\right),\quad
  \mu(c) =  F'(c)-\epsilon^2\frac{\partial^2 c}{\partial^2 x},\quad
  F(c) = c^2(1 - c)^2,
\end{equation}
where $c = c(x,t)$ is the unknown function,
$M = M(c) > 0$ is a diffusional mobility coefficient
(in all tests presented here we set~$M\equiv 1$),
$\mu = \mu(c,\partial{c}/\partial{x})$ is a chemical potential,
$\epsilon > 0$, $\epsilon = \text{const}$~is the diffuse interface width,
$t$~is the time, $x$~the space variable.

Equation~\eqref{eq:ch} is solved in a 1D domain~$\Omega = (0,1)$
for~$t \in (0,T]$.
Boundary conditions 
\begin{equation}\label{eq:bc}
  \frac{\partial c}{\partial x} = \frac{\partial \mu}{\partial x} = 0 ,
\end{equation}
are set on the domain boundary and, for ~$t=0$, 
initial condition is given,
\begin{equation}
\label{eq:ic}  
  c(x, 0) = c^0(x).
\end{equation}

From a physical point of view, equation~\eqref{eq:ch} can be derived as
follows~\cite{lee_2014}.
Consider a two-phase system with component concentrations~$c_{1,2}$, $c_1 + c_2 = 1$,
and assume that the free energy has a form, with~$c = c_1$,
\begin{equation}\label{eq:psi}
  \Psi[c,\nabla c] = \int\limits_\Omega\psi(c,\nabla c)\,dx,
\end{equation}
\begin{equation}\label{eq:psi:2}
  \psi(c,\nabla c) = F(c)+\frac{\epsilon^2}{2} |\nabla c|^2,\quad
  F(c) = c^2(1 - c)^2,
\end{equation}
where~$\psi$~is the system free energy density,
$F$~is the so-called double-well potential,
$\epsilon>0$~is a small parameter.
The first term in~\eqref{eq:psi:2} describes the ``dividing'' part 
of the free energy which ensures that the phases do not mix and,
hence, the regions with constant~$c_{1,2}$ emerge. 
The second term in~\eqref{eq:psi:2} allows to relate a given energy
to the diffuse interface.
Parameter~$\epsilon$ there defines the width of the diffuse
interface separating the ``pure'' phases.
Inclusion of the gradient terms makes the system energy depend not only 
on the quantity of each of the two components but also on the shape
of regions taken up by them.
Energy form~\eqref{eq:psi}, \eqref{eq:psi:2} arises 
in various models and is typical for weakly non-local (or gradient) thermomechanics.

It should be emphasized that potential~$F(c)$ in~\eqref{eq:psi:2} is empirical.
Its characteristic feature is the two minima at~$c=0$ and~$c=1$, corresponding
to the ``pure'' phases.
The state~$c=1/2$ corresponding to the maximum of~$F(c)$
is unstable and, hence, the mixing of the two phases is prevented
by the instability.
A more accurate and better physically motivated energy model
is the so-called logarithmic potential~\cite{lee_2014}.
It possesses essentially the same characteristic properties
as the energy form considered here, and the latter can be seen as
an approximation of the former.

\begin{figure}
  \centering
  \includegraphics[width = 0.48\textwidth]{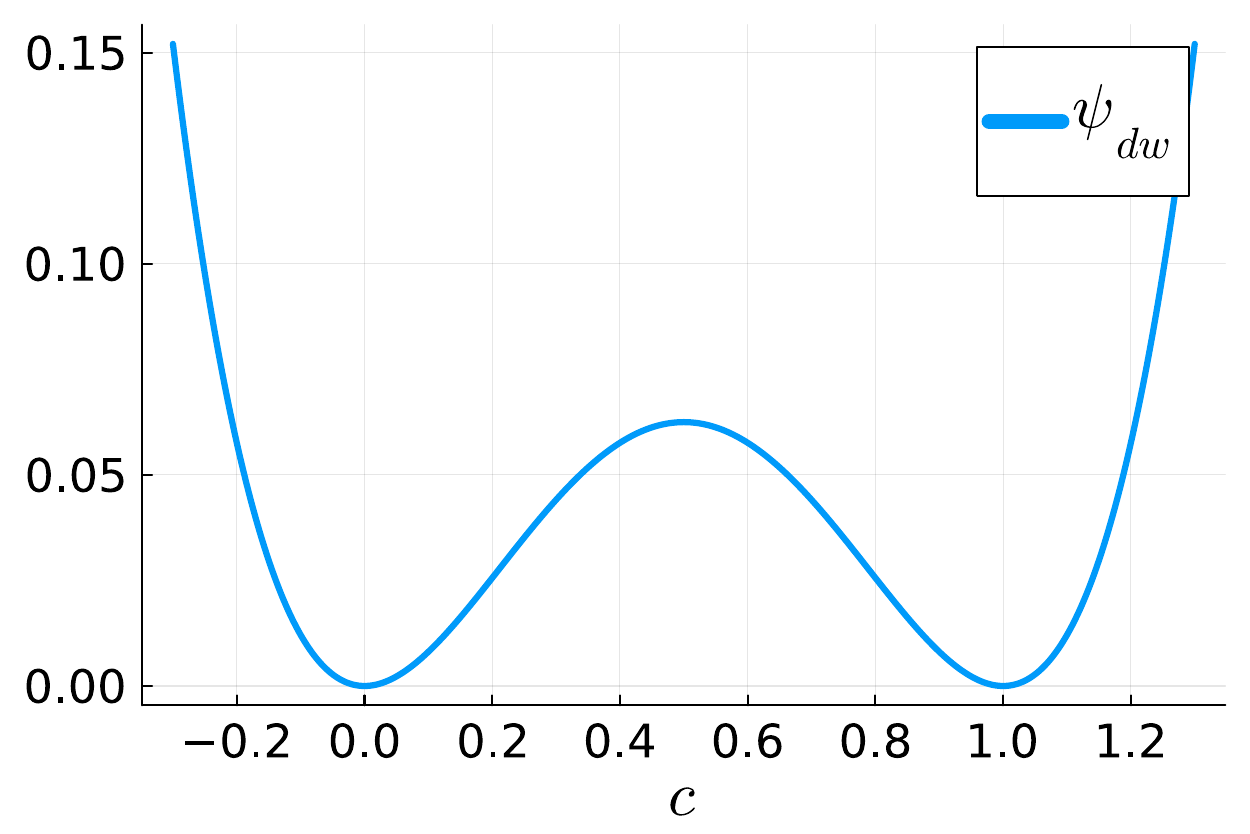}\hfill
  \includegraphics[width = 0.48\textwidth]{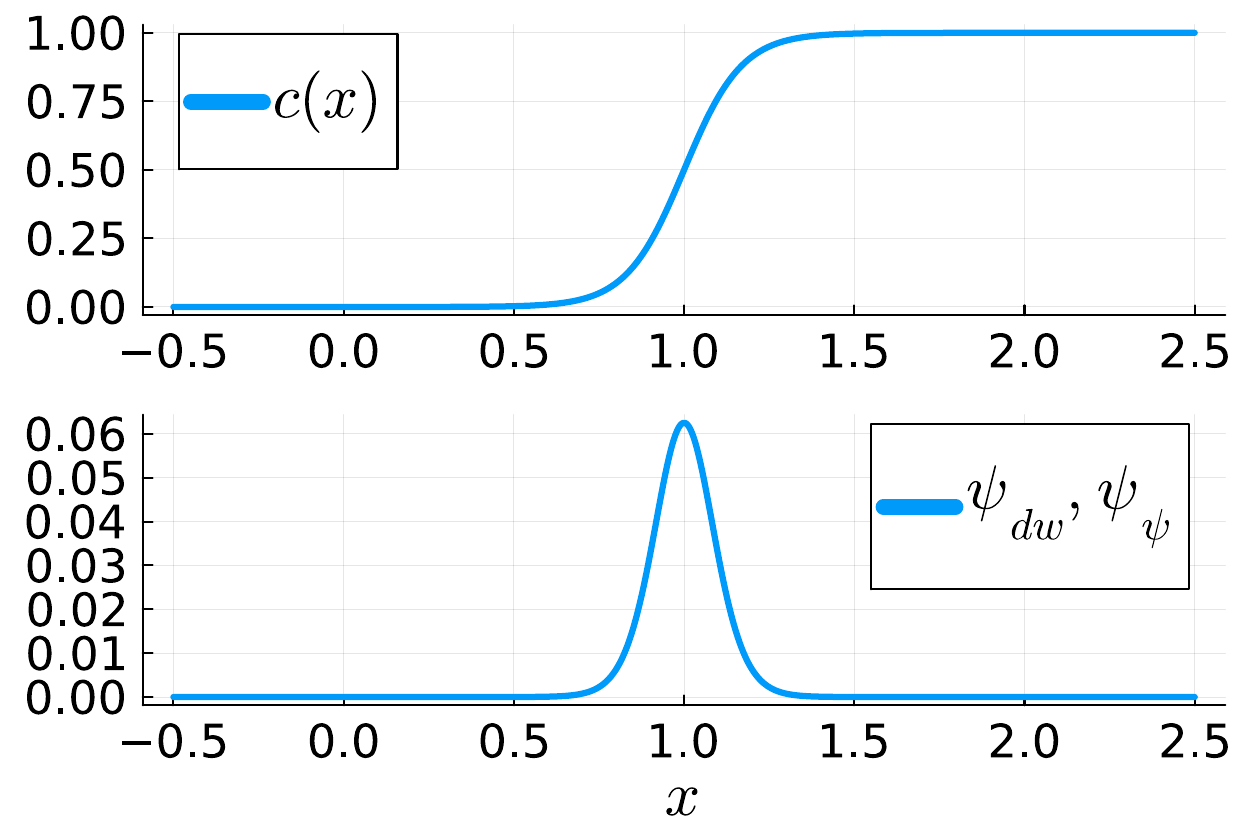}
  \caption{Left: an example of the ``separating'' free energy part~$\psi_{\text{dw}}$.
  Right: a typical solution and energy distribution.}
  \label{fig:psi}
\end{figure}

An illustration of these observations is given in Figure~\ref{fig:psi}.
In the left plot, the ``separating'' part~$\psi_\text{dw}\equiv F$ of the free energy 
is shown.
The two minima correspond to the ``pure'' phases,
i.e., to the states~$c=c_1=0$, $c_2=1$ and~$c=c_1=1$, $c_2=0$.
The state in the neighborhood of~$c=1/2$ is unstable.
In the right plot of the figure, a typical stationary solution of the system
and energy distribution are shown.
As can be seen, the gradient part of the free energy is nonzero 
only in the diffuse interface regions.
At the same time, the separating part of the free energy is zero outside
the diffuse interface regions.

Furthermore, one can show that if~$c$ is a conservative quantity (as is
in our case), the kinetic equation describing the field evolution~$c = c(x,t)$
reads~\cite{hohenberg_1977,lee_2014}:
\begin{equation}\label{eq:ch:2}
  \frac{\partial c}{\partial t}
  = -\nabla\mathbf{\cdot}\left( -M \nabla \mu \right), \quad
  \mu =\frac{\delta\psi[c,\nabla c]}{\delta c},
\end{equation}
where~$\delta(\cdot)/\delta c$ is the functional derivative
(Gateaux derivative).
For~$M$ constant, due to~\eqref{eq:psi:2}, the last relation takes form
\[
  \frac{1}{M}
  \frac{\partial c}{\partial t}
  = \Delta\mu, \quad
  \mu = -\epsilon^2\Delta{c} + F'(c),
\]
where $\Delta$ is the Laplacian.
In the 1D case, under consideration here, this coincides with~\eqref{eq:ch}.

As can also be shown~\cite{demello_2004,lee_2014}, equation~\eqref{eq:ch}
has solutions $c = c(x,t)$ with the following fundamental properties:
\[
  \frac{d}{dt} \int\limits_{\Omega} c = 0,
  \qquad
  \frac{d}{dt} \Psi[c,\nabla c] \leqslant 0,
\]
which means that $c$ is conserved and 
the free energy is a non-increasing function of
$c$ being the solution of~\eqref{eq:ch} or~\eqref{eq:ch:2}
and satisfying homogeneous boundary conditions~\eqref{eq:bc}.

These properties of the analytical solution immediately 
put restrictions on a numerical scheme for solving~\eqref{eq:ch}:
it should be conservative (which is relatively easy to achieve)
and
energy-stable, i.e., delivering a discrete analogue of the 
last inequality (which is more difficult to achieve).
In addition, the discrete problem should be uniquely solvable for 
a given time step size.  In practice, this may turn out to be a rather 
nontrivial task because the potential~$F(c)$ is not a convex function
and, hence, a discrete solution may be found in a local minimum
which is not global and not achievable from an initial state.

\section{Numerical algorithms}\label{sec:algo}
This section deals with numerical algorithms applied to solve 
the Cahn-Hilliard equation.  In the first subsection
classical time integration schemes are considered, having various degree of 
implicitness.  
The second subsection describes the Eyre splitting for
the time integration of the Cahn-Hilliard equations. 
New algorithms based on the LIM scheme are 
presented in the third subsection.

\subsection{Classical time integration schemes}
\label{sec:class_alg}
We now consider a number of time integration schemes taking part
in the numerical comparisons presented below.
All these schemes are well known and described in numerous literature
and, in particular, in~\cite{lee_2013}.  
All schemes considered in this subsection approximate
equation~\eqref{eq:ch} and are conservative.  

Throughout this paper space derivatives are discretized by finite differences.
Assume that a uniform grid with grid size~$h=1/N$
and nodes~$x_i = (i-1/2) h$, $i=\overline{0,N+1}$,
with $N$~being the number of grid cells, is used to discretize
the equation in the domain~$\Omega$.
The domain boundaries~$x=0$ and~$x=1$ correspond to points
$x_{1/2} = (x_0+x_1)/2=0$ and~$x_{N+1/2} = (x_{N}+x_{N+1})/2 = N h$
and the grid nodes~$x_1,x_2,\ldots,x_N$~lie within the domain.
For time discretization a uniform grid with time step size~$\tau$ 
and nodes~$t_n = n \tau $, $n=0,1,2,\ldots$, is used.
Thus, the solution of the discretized problem is defined in
points~$(x_i,t_j)$,
$i=\overline{1,N}$, $j=0,1,\ldots$.
The corresponding values of the grid functions~$c_h$ are denoted
by~$c_i^n = c(x_i,t_n)$. 

Let~$\Delta_h$ be a standard three-point finite-difference 
approximation of the Laplacian. 
In the 1D case under consideration, its values in the nodes~$i=\overline{1,N-1}$ 
of the space grid are 
\[
  \Delta_h c_i = \frac{1}{h^2} \left( c_{i+1} - 2 c_i + c_{i-1} \right).
\]
The finite difference equations 
\begin{equation}\label{eq:bc:2}
  c_0 - c_1 = 0,\quad c_N - c_{N-1} = 0
\end{equation}
approximating the homogeneous Neumann boundary conditions~\eqref{eq:bc}
are related to the boundary nodes~$i=0,N$.
Then the resulting discrete Laplacian, with incorporated 
discretized boundary conditions, is denoted by~$\Delta_h$.

\textbf{Explicit Euler (EE)} scheme, the simplest 
discretization of~\eqref{eq:ch}, reads
\begin{gather*}
  \frac{1}{M}\frac{c^{n+1}_{i} - c^{n}_{i}}{\tau} = \Delta_h\mu^{n}_{i},\quad
\mu^{n}_{i} = F'(c^{n}_{i}) - \epsilon^2\Delta_h c^{n}_{i}.
\end{gather*}
The scheme is first order accurate in time and second order
accurate in space, is stable for~$\tau \sim h^4$
and is not energy stable.

\textbf{Implicit Euler (IE)} scheme reads
\begin{gather*}
  \frac{1}{M}\frac{c^{n+1}_{i} - c^{n}_{i}}{\tau} = \Delta_h\mu^{n+1}_i,\quad
  \mu^{n+1}_{i} = F'(c^{n+1}_{i}) - \epsilon^2\Delta_h c^{n+1}_{i}.
\end{gather*}
The scheme is first order accurate in time and second order
accurate in space and is nonlinear 
(in the sense that a nonlinear equation has to be solved at each
time step to determine~$c^{n+1}_{i}$).
It is conditionally (i.e., for a sufficiently small~$\tau$)
gradient stable and conditionally uniquely solvable, see
Table~\ref{tab:schemes}.

\textbf{Crank-Nicolson (CN)} scheme can be written as
\begin{gather*}
\frac{1}{M}\frac{c^{n+1}_{i} - c^{n}_{i}}{\tau} = \Delta_h \mu^{n+1/2}_{i},\quad
\mu^{n+1/2}_{i} = \frac{1}{2}(\mu^{n}_{i}+\mu^{n+1}_{i}),
\end{gather*}
with~$\mu_i^{n}$, $\mu_i^{n+1}$ defined above.
The scheme is second order accurate in time and second order
accurate in space and is nonlinear.
It is conditionally 
gradient stable and conditionally uniquely solvable, see
Table~\ref{tab:schemes}.

\textbf{Semi-implicit Euler (SIE)} reads
\begin{gather*}
  \frac{1}{M}\frac{c^{n+1}_{i} - c^{n}_{i}}{\tau} = \Delta_h\mu^{n+1}_{i},\quad
  \mu^{n+1}_{i} = F'(c^{n}_{i}) - \epsilon^2 \Delta_h c^{n+1}_{i}.
\end{gather*}
The scheme is first order accurate in time and second order
accurate in space and is linear.
It is unconditionally uniquely solvable and is not energy stable, see
Table~\ref{tab:schemes}.
%

\subsection{Eyre splitting}
\label{sec:Eyre}
The time integration schemes considered up to now 
are not energy stable and/or not uniquely solvable for
arbitrary large time step sizes~\cite{demello_2004,vollmayr_2003}.
An efficient and general way to design time integration schemes,
which are both uniquely solvable and energy stable for arbitrary 
time step sizes, is proposed in papers of David J. Eyre~\cite{eyre_1997,eyre_1998}.
Since we essentially employ this approach to construct new schemes,
we discuss it in more details here.

The essence of the Eyre approach is as follows.
The system free energy density~$\psi$, c.f.~\eqref{eq:psi},
is not a convex function.   This is because the dividing part of the 
free energy is a double-well potential.
Assume that~$\psi$ can be split as
\begin{equation}\label{eq:split}
  \psi = \psi_\text{c} + \psi_\text{e},
\end{equation}
where~$\psi_{\text{c}}$ and $-\psi_\text{e}$~are \emph{convex} functions.
The subindices~``c'' (meaning ``contraction'') and~``e'' (meaning ``expansion'')
reflect their physical meaning. 
Here the convexity of a function is understood as a positive semi-definiteness
of its Hessian, see~\cite{vollmayr_2003},
i.e., a function~$f(x)$ is convex provided that~$f''(x)\equiv d^2f/dx^2\geqslant 0$.
Then, the semi-discrete (discrete in time) approximations 
\begin{gather}\label{eq:eyre}
\frac{1}{M}  \frac{c^{n+1} - c^n}{\tau}
  = \Delta\left[
    \mu_\text{c}(c^{n+1}) + \mu_\text{e}(c^{n})
    \right],
\end{gather}
with
\[
  \mu_\text{c,e} = \frac{\delta \psi_\text{c,e}}{\delta c}, \quad
  \mu = \mu_\text{c} + \mu_\text{e},
\]
are gradient stable and 
\[
  \psi(c^{n+1}) \leqslant \psi(c^n)
\]
holds for any time step size~$\tau> 0$. In addition, the convexity
of the energy part$\psi_\text{c}$ corresponding to the implicitly treated term,
ensures the unique solvability of~\eqref{eq:eyre}.

Consider the specific form of the energy functional~\eqref{eq:psi}.
Taking into account that the gradient term in~\eqref{eq:psi}
belongs to the convex part of the energy,
we can take in~\eqref{eq:split}
\[
  \psi_\text{c} = \frac{1}{2}|\nabla c|^2 + F_\text{c},\quad
  \psi_\text{e} = F_\text{e}, \quad F_\text{c} + F_\text{e} = F.
\]
Thus, constructing a suitable splitting of the system energy is reduced
to splitting the dividing part of the energy~$F$ in such a way that
\begin{equation}\label{eq:split:1}
  F = F_\text{c} + F_\text{e}, \quad 
F''_\text{c} \geqslant 0,
\quad
-F''_\text{e} \geqslant 0.
\end{equation}
Then, the chemical potential takes form
\[
  \mu = \mu_\text{c} + \mu_\text{e},\quad
  \mu_\text{c} = -\epsilon^2\Delta c + F'_\text{c},\quad \mu_\text{e} = F'_\text{e},
\]
whereas the semi-discrete system~\eqref{eq:eyre} can be written as
\begin{equation}\label{eq:eyre:2}
  \frac{1}{M}\frac{c^{n+1} - c^n}{\tau}
  = -\epsilon^2\Delta^2 c^{n+1}
  + \Delta F'_\text{c}(c^{n+1}) + \Delta F'_\text{e} (c^n).
\end{equation}

Splitting of the form~\eqref{eq:split:1} can be constructed in a number of ways.
One of them is to write the dividing part~$F$ of the free energy 
as a sum of two terms according to~\eqref{eq:split}.
Note that the gradient part of the free energy~\eqref{eq:psi:2} 
should be handled together with the convex part~$\psi_\text{c}$.
Another way is to choose a regularizing function~$F_\text{r}$ such that
\[
  F_\text{c} = F + F_\text{r},\quad F_\text{e} = F - F_\text{r}.
\]
Although both ways are essentially quite similar,
it is convenient to distinguish them due to various possible
generalizations.

Consider the specific form~\eqref{eq:psi:2} of~$F$
used in this work.  We have
\[
  F(c) = c^2 (1-c^2), \quad
  F'(c) = 2 c (2 c^2 - 3 c + 1),\quad
  F''(c) = 2 (6 c^2 - 6 c + 1).
\]
Equation $F''(c)=0$ has roots~$c_{1,2} = \left(3\pm\sqrt{3}\right)/6$,
and, hence, $F''(c) < 0$
for~$\left(3-\sqrt{3}\right)/6 < c < \left(3+\sqrt{3}\right)/6$\,;
$F''(c) \geqslant 0$
for~$c\leqslant\left(3-\sqrt{3}\right)/6$
and~$c\geqslant\left(3+\sqrt{3}\right)/6$.
Function $F''(c)$ attains its minimum, equal to~$-1$, at~$c=1/2$.
The maxima of~$F''(c)$ for $c\in[0,1]$ are attained at~$c=0$ and~$c=1$
and equal~$2$.
Thus, to ensure that~$F_\text{c}$ and~$-F_\text{e}$ are convex,
it suffices to set
\[
  F_\text{c} = F(c) + \frac{1}{2}c^2,\quad F_\text{e} = -\frac{1}{2}c^2.
\]
This leads to \textbf{non-linearly stabilized splitting scheme (NLSS)},
\begin{gather}\label{eq:nlss}
  \frac{1}{M} \frac{c^{n+1}_{i} - c^{n}_{i}}{\tau}
  = \Delta_h \mu^{n+1},\quad
  \mu^{n+1}_{i} = \left[F'(c^{n+1}_{i})  + c^{n+1}_{i}\right] - c^{n}_{i}  -\epsilon^2\Delta_h c^{n+1}_{i}.
\end{gather}
The scheme can also be written as
\[
  \frac{1}{M}\frac{c^{n+1}_{i} - c^{n}_{i}}{\tau}
  = \left(\Delta_h - \epsilon^2\Delta_h^2\right) c^{n+1}_i + \Delta_h F'(c^{n+1}_i)
    - \Delta c^n_i.
\]
It is first order accurate in time, second order in space,
implicit, nonlinear in~$c^{n+1}$ and energy stable for 
arbitrary~$\tau>0$.

Another option is to set
\[
  F_\text{c} = c^2,\quad F_\text{e} = -c^2 + F(c),
\]
which leads to \textbf{linearly stabilized splitting scheme (LSS)},
\begin{gather*}
  \frac{1}{M}
  \frac{c^{n+1}_{i} - c^{n}_{i}}{\tau}
  = \Delta_h \mu^{n+1}_{i},\quad
  \mu^{n+1}_{i} = 2c^{n+1}_{i} + \left[F'(c^{n}_{i}) -
    2c^{n}_{i}\right] - \epsilon^2\Delta_h c^{n+1}_{i}.
\end{gather*}
It may be convenient to write it as
\[
  \frac{1}{M}\frac{c^{n+1}_{i} - c^{n}_{i}}{\tau}
  = \left(2\Delta_h - \epsilon^2\Delta_h^2\right) c^{n+1}_i + \Delta_h F'(c^{n}_i)
    - 2\Delta c^n_i.
\]
The scheme is first order accurate in time, second order in space, 
implicit and, provided that splitting parameters are chosen properly~\cite{wise_2009},
energy stable for arbitrary~$\tau>0$.
Unlike the NLSS scheme~\eqref{eq:nlss},
it is linear with respect to~$c^{n+1}$.
A three-parameter family of linearly stabilized schemes,
including the LSS scheme, is analyzed in~\cite{vollmayr_2003}.

Properties of considered time integration schemes are summarized in
Table~\ref{tab:schemes}.
\begin{table}[t]
\caption{Properties of time integration schemes according to~\cite{demello_2004,vollmayr_2003,eyre_1997,eyre_1998,wise_2009}
    ($^*$~provided parameters are chosen properly, see~\cite{wise_2009}).}
\label{tab:schemes}
\addtolength{\tabcolsep}{0.4em}
  \centering
  \bgroup
  \def\arraystretch{1.5}
\begin{tabular}{ l l l l}
\hline\hline
Scheme              & Linearity  & Gradient stability   & Solvability                                 \\ \hline
Explicit Euler (EE) & Yes        & No                   & Yes \\ \hline
Implicit Euler (IE) & No         & Conditional, $\tau\leqslant \frac{1}{4}h^2$ 
                                                        & $\tau \leqslant \frac{1}{18}h^2$ \\ \hline
Crank-Nicolson (CN) & No         & Conditional          & $\tau \leqslant \frac{1}{9}h^2$  \\ \hline
Semi-implicit Euler (SIE)        
                    & Yes        & No                   & Yes  \\ \hline
Nonlinearly stabilized splitting (NLSS)
                    & No         & Unconditional        & Yes  \\ \hline
Linearly stabilized splitting (LSS)
                    & Yes        & Unconditional$^*$    & Yes   \\  \hline
  \end{tabular}
  \egroup
\addtolength{\tabcolsep}{-1em}
\end{table}

For a fully discrete scheme (i.e., when discretization in time and in space
is applied), its gradient stability implies that the free energy 
functional~\eqref{eq:psi} is approximated in a certain way.
In this work such an approximation, determined by the space disretization,
reads
\begin{equation}
\label{psi_h}  
  \psi_h(\cc)
  = h \sum\limits_{i=1}^N F'(c_i)
  + \frac{1}{2}\epsilon^2 h \sum\limits_{i=1}^{N} \frac{(c_{i+1}-c_i)^2}{h^2}.
\end{equation}
%

Let a~$N\times N$ matrix~$\AB$ be the matrix of the discussed
above discretization $-\Delta_h$ of the operator~$-\Delta$,
with incorporated homogeneous Neumann boundary conditions~\eqref{eq:bc:2}.
The space discretization reduces the initial boundary-value 
problem~\eqref{eq:ch}--\eqref{eq:ic} to an initial value 
problem
\begin{equation}
\label{eq:ivp}
\frac 1 M \cc' = 
-\AB\left(F'(\cc) + \epsilon^2\AB\cc\right),
\quad \cc(0)=\cc^0,
\end{equation}
where the entries of the vector function $\cc(t):\Rr\rightarrow \Rr^N$ are 
the sought after solution values on the space grid and
the vector $\cc^0$ contains the grid values of the given function $c^0(x)$.

Let $\cc^n$ be an~$N$ vector containing the numerical 
solution values for time level~$n$.
Then, the time integration schemes discussed above can be written
in a compact form as follows.
\begin{itemize}
\item Explicit Euler (EE):
  \begin{equation}\label{explicit}
    \frac{\cc^{n+1} - \cc^{n}}{\tau} = -\AB\left(F'(\cc^{n})+\epsilon^2\AB\cc^{n}\right);
  \end{equation}

\item implicit Euler (IE):
  \begin{equation}\label{implicit}
    \frac{\cc^{n+1} - \cc^{n}}{\tau} = -\AB\left(F'(\cc^{n+1})+\epsilon^2\AB\cc^{n+1}\right);
  \end{equation}

\item Crank-Nicolson (CN):
  \begin{equation}\label{CN}
    \frac{\cc^{n+1} - \cc^{n}}{\tau} = -\frac{1}{2}\AB\left(F'(\cc^{n+1})+\epsilon^2\AB\cc^{n+1}+F'(\cc^{n})+\epsilon^2\AB\cc^{n}\right);
  \end{equation}
  
\item Semi-implicit Euler (SIE):
  \begin{equation}\label{SIE}
    \frac{\cc^{n+1} - \cc^{n}}{\tau} = -\AB\left(F'(\cc^{n})+\epsilon^2\AB\cc^{n+1}\right);
  \end{equation}

\item nonlinearly stabilized splitting (NLSS)
\begin{equation}\label{NLSS}
  \frac{\cc^{n+1} - \cc^{n}}{\tau} = -\AB\left(F'(\cc^{n+1})+\cc^{n+1}-\cc^{n}+\epsilon^2\AB\cc^{n+1}\right);
\end{equation}

\item linearly stabilized splitting (LSS)
\begin{equation}\label{LSS}
    \frac{\cc^{n+1} - \cc^{n}}{\tau} = -\AB\left(F'(\cc^{n})- 2\cc^{n}+2\cc^{n+1}+\epsilon^2\AB\cc^{n+1}\right).
  \end{equation}
  
\end{itemize}

We now briefly discuss computational complexity of these schemes.
In implicit nonlinear schemes (i.e., IE, CN, NLSS), 
each time step a nonlinear system of equations has to be solved, 
e.g., by a Newton iteration.
In turn, each Newton iteration involves solving the Jacobian linear system.
Implicit linear schemes SIE and LSS require solving a single
linear system per time step. 
Explicit schemes involve neither nonlinear nor linear system solution
and have lowest costs per time step.
Although implicit time integration schemes allow to use significantly
larger time step sizes than explicit schemes,
solving large linear systems, within Newton iterations or independently,
may become a computationally expensive task.
This is because efficiency in sparse direct solvers, as well as in
preconditioned iterative solvers, is difficult to retain 
on modern high-performance hybrid CPU/GPU platforms. 

The aim of this work is to develop an efficient numerical algorithm
which, on one hand, is gradient stable and allows large time step
sizes and, on the other hand, can be efficiently implemented
on modern high-performance hybrid platforms. 
Our approach is based on the combination of the Eyre splitting,
in particular the LSS scheme, with 
the local iteration modified (LIM) scheme.
Then, the Eyre splitting allows to use large time step sizes
while keeping gradient stability, whereas 
the LIM mechanism provides a computationally efficient 
alternative to solving linear systems in implicit schemes.

\subsection{Local iteration modified (LIM) scheme}
The local iteration modified (LIM) scheme can be viewed as a special
explicit scheme where a Chebyshev polynomial of 
a sufficiently high order~$p$ is employed to ensure
stability and monotonicity for a given time step size $\tau>0$.
We now describe a LIM variant based on the linearly stabilized splitting 
(LSS) scheme~\eqref{LSS}, which we write as
\begin{equation}
\label{LSS2}
\cc^{n+1} = \left[ \IB + \tau \widehat{\AB} \right]^{-1}\widehat{\ff}^n,
\quad \widehat{\AB} = \AB (2\IB + \epsilon^2 \AB),
\quad
\widehat{\ff}^n = \cc^n + \tau \AB (2\cc^n - F'(\cc^n)),
\end{equation}
where $\IB$ is the $N\times N$ identity matrix.
We emphasize that the inverse matrix here is, of course,
never computed, rather a linear system with this matrix is solved.
Within the LIM approach the inverse matrix is replaced by a specially chosen
Chebyshev polynomial.
Let $\lambda_{\infty}$ be an upper bound for the largest eigenvalue of~$\widehat{\AB}$
(in practice one can set $\lambda_{\infty}:=\|\widehat{\AB}\|_1=\max_j\sum_i|\widehat{a}_{ij}|$).
The polynomial order~$p$ for which stability takes place is chosen as 
\begin{equation}
\label{p}
p = \left\lceil \frac\pi4\sqrt{\tau\lambda_{\infty} + 1}\,\right\rceil,
\end{equation}
where $\lceil x \rceil$, for $x\in\Rr$, denotes the smallest integer 
greater than or equal to~$x$.  
Let the Chebyshev polynomial roots 
$$
\{ \beta_m, \; m=1,\ldots,p \} = \left\{ \cos\pi\frac{2i-1}{2p}, \; i=1,\dots,p\right\}
$$
be ordered in such a way that no numerical instability occurs 
for arbitrarily large~$p$
and $\beta_1=\cos({\pi}/{2p})$ is the first root.
We set $z_1=\beta_1$ and define Chebyshev polynomial parameters
$$
a_m=\frac{\lambda_{\infty}}{1+z_1}(z_1 - \beta_m), \quad m=1,\dots,p.
$$
Denoting $\yy^{(0)}=\cc^n$, we obtain the LIM solution~$\cc^{n+1}$ 
on the next time level by carrying out $2p-1$ Chebyshev iterations as
\begin{equation}
\label{LI-M}
\begin{aligned}
\yy^{(m)} &= \frac1{1+\tau a_m}\left(
\cc^n + \tau a_m\yy^{(m-1)} + \tau(\widehat{\ff}^{n} - \widehat{\AB}\yy^{(m-1)})
\right),
\quad m=1,\dots,p,
\\
\yy^{(p+m-1)} &= \frac1{1+\tau a_m}\left(
\cc^n + \tau a_m\yy^{(m-1)} + \tau(\widehat{\ff}^{n} - \widehat{\AB}\yy^{(m-1)})
\right),
\quad m=2,\dots,p,
\\
\cc^{n+1} &= \yy^{(2p-1)}.
\end{aligned}
\end{equation}
Note that $a_1=0$ and it is easy to check that the first iteration solution~$\yy^{(1)}$
is the explicit Euler solution on the next time level.
In this sense, the Chebyshev iterations actually start
at the second step by computing~$\yy^{(2)}$.
If only the first sweep of iterations is carried out in~\eqref{LI-M},
with $\yy^{(1)}$, \dots, $\yy^{(p)}$ computed, then setting $\cc^{n+1} = \yy^{(p)}$
we get a time integration scheme which is known as the regular (non-modified) 
local iteration scheme (the LI scheme).  As already mentioned,
in the regular LI scheme stability of the numerical solution
is guaranteed for parabolic problems but the solution nonnegativity 
is not.
As one can see, in the LIM scheme~\eqref{LI-M} the second Chebyshev iteration
sweep, with $\yy^{(p+1)}$, \dots, $\yy^{(2p-1)}$ being computed,
the same parameters~$a_m$ are employed, except that the first iteration 
with~$a_1=0$ (which gives the explicit Euler solution) is not carried out.
One can show~\cite{zhukov_2011} that formulas~\eqref{LI-M} 
can be written in an operator form as
\begin{equation}
\label{LI-Moper}
\cc^{n+1} = (\IB-\FB_p^2)\left[ \IB + \tau \widehat{\AB} \right]^{-1}
\widehat{\ff}^n,
\end{equation}
where $\widehat{\AB}$ and $\widehat{\ff}^n$ are defined in~\eqref{LSS2}
and $\FB_p$ is the Chebyshev polynomial operator,
$$
\FB_p = \prod_{m=p}^{m=1}\left(
\IB-\frac1{1+\tau a_m}(\IB + \tau \widehat{\AB})
\right).
$$
By replacing $\FB_p^2$ in~\eqref{LI-Moper} by $\FB_p$ we obtain
an operator form of the (regular) LI scheme. 
Thus, as we see, squaring the Chebyshev polynomial leads to
monotonicity in the LIM scheme.

Since scheme~\eqref{LI-M} under consideration is derived from
the linearly stabilized splitting~\eqref{LSS}, we call
this scheme LIM-LSS.  
To discover the Eyre splitting effect in this scheme,
we include another scheme the numerical experiments
presented below.  It is obtained by linearizing implicit Euler
scheme~\eqref{implicit}.
We linearize the scheme by approximating the nonlinear implicit
scheme in~\eqref{implicit} as
\begin{equation}
\label{lin}
F'(\cc^{n+1}) \approx F'(\cc^n) + J_n (\cc^{n+1} - \cc^n),
\end{equation}
where $J_n$ is the Jacobian of the mapping $F'$ evaluated at $\cc^n$.
Substituting approximation~\eqref{lin} in~\eqref{implicit} leads to
a linearized implicit Euler (LIE) scheme
\begin{equation}
\label{lin_impl}
\cc^{n+1} = \left[ \IB + \tau \widetilde{\AB}_n \right]^{-1}\widetilde{\ff}^n,
\quad \widetilde{\AB} = \AB (J_n + \epsilon^2 \AB),
\quad
\widetilde{\ff}^n = \cc^n + \tau \AB (J_n\cc^n - F'(\cc^n)),  
\end{equation}
where, just as in scheme~\eqref{LSS2}, a linear system solution
is understood, so that the inverse matrix is not computed.
Comparing schemes~\eqref{LSS2} and~\eqref{lin_impl}, one can notice
that replacing $J_n$ in~\eqref{lin_impl} by $2\IB$ leads to scheme~\eqref{LSS2}.

The LIM scheme based on~\eqref{lin_impl} can be obtained in exactly
the same as the LIM scheme based on~\eqref{LSS2}.  
One can simply repeat the derivations given above, replacing 
$\widehat{\AB}$ and $\widehat{\ff}_n$ by 
$\widetilde{\AB}_n$ and $\widetilde{\ff}_n$, respectively.  
We call the LIM scheme derived in this way 
the LIM-LIE  (LIM linearized implicit Euler) scheme.

A computational efficiency of the LIM schemes can be established,
based on the number of required iterations~\eqref{p}
and the following reasoning~\cite{LokLok1986, zhukov_2011}.
If the time interval for which the problem has to be solved
is increased by a factor~$s$ then computational costs
of an explicit scheme should also grow by a factor~$s$
(as the number of required time steps is $s$~times larger). 
The same increase in computational costs is also observed
if the upper spectral bound $\lambda_{\infty}$ is increased
by a factor~$s$ (since a time step size in an explicit scheme
is usually bound to a stability condition of the form 
$\tau\leqslant 2/\lambda_{\infty}$).  
Recall that for the Cahn-Hilliard equation, in general,
we have $\lambda_{\infty}\sim h^{-4}$.
In the LIM schemes the costs grow slower:
an increase of the time interval or the upper spectral bound $\lambda_{\infty}$
by a factor~$s$ leads to an increase in computational work
approximately by a factor $\sqrt{s}$
(as the number of required Chebyshev iterations 
$p\sim (\tau\lambda_{\infty})^{1/2}$, see~\eqref{p}).
These estimates make clear that the gain in computational
costs provided by the LIM schemes increases with the problem
size $N=1/h$~\cite{LokLok1986, zhukov_2011}.
Nevertheless, it should also be taken into account
that the gain actually observed in practice can be
restricted by other factors, such as accuracy requirements.

\section{Computational experiments}\label{sec:num}

\subsection{Gradient stability tests}
In the tests presented here initial-value problem~\eqref{eq:ch}--\eqref{eq:ic}
is solved for the 1D Cahn-Hilliard equation.
Space discretization is done by the standard second order finite
differences on a uniform grid,
as discussed in the beginning of Section~\ref{sec:class_alg}.
This reduces the original initial-boundary
value problem to an initial value problem~\eqref{eq:ivp}
to be solved by the time integration schemes~\eqref{explicit}--\eqref{NLSS}
and new local iteration schemes~\eqref{LI-M} 
based on~\eqref{LSS} and~\eqref{lin_impl} (the LIM-LSS and LIM-LIE schemes).
The parameter $\epsilon$ is set to one of the following two values:
\begin{align}
\label{eps_h}
\epsilon &= \epsilon_4(h), \quad
\epsilon_m(h) \equiv \frac{h m}{2\sqrt{2} \, \mathrm{arth} (9/10)},
\\
\label{eps_fix}
\epsilon &= \epsilon_4(1/64),
\end{align}
i.e., $\epsilon$ either depends on the space grid size~$h$
(formula~\eqref{eps_h},
or set to a fixed value $\epsilon=\epsilon_4(1/64)$
of the grid~$N=64$ (formula~\eqref{eps_fix}).
Note that $m=4$ in formula~\eqref{eps_h} defines a characteristic
number of the grid cells on which a stationary solution changes
from its minimum to maximum values.

For all the tested time integration schemes, largest time step sizes
for which the gradient stability is observed are reported in
Table~\ref{t:gr_st}.
In these runs, the grid-dependent $\epsilon$ values~\eqref{eps_h}
are used and a scheme is considered to be gradient stable
provided that the discrete energy~\eqref{psi_h} does not increase
more than 1\% each time step, i.e., for each time step~$n$ holds
\begin{equation}
\label{psi_hn}
\psi_h(\cc^{n+1})\leqslant 1.01\psi_h(\cc^n).
\end{equation}
Gradient stability is tested on several initial value vectors~$\cc^0$,
where each vector component is taken to be
an independent and identically distributed random variable
in~$[0,1]$, rounded up to the second digit after the decimal point. 
Table~\ref{t:gr_st} also contains results for the
linearized implicit Euler (LIE) scheme~\eqref{lin_impl}. 

As can be seen in the table, only Eyre splitting based schemes
turn out to be unconditionally gradient stable.
In the tests, a scheme is considered to be unconditionally gradient stable
if gradient stability is observed for $\tau \leqslant 500$.
Furthermore, we note that maximum time step size in the
implicit scheme exceeds the maximum time step size in the
explicit schemes approximately by the same factor ($\approx 60$)
for all~$h$.
Comparing the LIE and IE schemes, we see that linearization does
not corrupt the gradient stability of the implicit scheme.
The LIM-LSS scheme conserves the unconditional gradient stability of the
LSS scheme, and the LIM-LIE scheme inherits the conditional
gradient stability of the LIE scheme.
One of the aims of this work is to show that implicit schemes
can be successfully replaced by stable explicit schemes.
Having this in mind, we may conclude from Table~\ref{t:gr_st} that
both LIM schemes work well.
This, of course, does imply that increasing the time step
size~$\tau$ in both LIM-LSS and LIM-LIE schemes leads to an increase
of the Chebyshev iterations required per time step.
We also note that the values presented in our Table~\ref{t:gr_st}
are close to the values reported in Table~1 from paper~\cite{lee_2013}.

\begin{table}
\caption{Largest time step size $\tau$ for which
time integration schemes are gradient stable,
i.e., condition~\eqref{psi_hn} holds for $\epsilon=\epsilon_4(h)$}
\label{t:gr_st}
\begin{center}
\begin{tabular*}{0.8\textwidth}{@{\extracolsep{\fill} } lcccc}
\hline\hline
$h$ & $1/32$ & $1/64$ & $1/128$ & $1/256$   \\
\hline
EE            & $8.8\times10^{-5}$    & $2.1\times10^{-5}$     & $5.4\times10^{-6}$     & $1.3\times10^{-6}$     \\
IE            & $5.1\times10^{-3}$    & $1.3\times10^{-3}$     & $3.0\times10^{-4}$     & $8.1\times10^{-5}$     \\
CN            & $2.4\times10^{-3}$    & $6.8\times10^{-4}$     & $1.7\times10^{-4}$     & $3.8\times10^{-5}$     \\
SIE           & $2.2\times10^{-3}$    & $5.7\times10^{-4}$     & $1.9\times10^{-4}$     & $4.0\times10^{-5}$      \\
LSS           & $\infty$    & $\infty$     & $\infty$    & $\infty$      \\
NLSS          & $\infty$    & $\infty$     & $\infty$    & $\infty$      \\
LIM-LSS       & $\infty$    & $\infty$     & $\infty$    & $\infty$      \\
LIM-LIE       & $9.9\times10^{-3}$    &  $2.8\times10^{-3}$  & $5.8\times10^{-4}$ & $1.2\times10^{-4}$ \\
LIE           & $5.8\times10^{-3}$    &  $1.6\times10^{-3}$  & $3.2\times10^{-4}$ & $8.3\times10^{-5}$ \\
\hline
\end{tabular*}
\end{center}
\end{table}

\subsection{Accuracy and efficiency tests}
The tests discussed here are aimed to evaluate accuracy and
computational efficiency of the considered schemes and,
in particular, to discover whether the LIM schemes provide
a gain in computational costs with respect to the explicit Euler
(EE) scheme.
Since the costs per time step in the EE~scheme are minimal,
it is clear that other time integration schemes can deliver
a higher computational efficiency only if the time step size is
increased.
However, increasing the time step size is possible only
provided the delivered numerical accuracy remains
within the allowed tolerance.
We estimate accuracy of a time integration scheme
by comparing its solution at the final time $t=T$
to the reference solution $\cc_{\text{ref}}(T)$,
which is computed for each space grid with a tiny time
step size $\tau = 10^{-9}$.
The initial value vector is set to have random entries,
in the same way as in the previous tests, then fixed and
used for a current space grid in all time integration schemes.
For each scheme the achieved accuracy is measured as
a relative error norm
\begin{equation}
\label{err} 
\frac{\|\cc^n-\cc_{\text{ref}}(T)\|}{\|\cc_{\text{ref}}(T)\|},
\quad n= n_{\text{final}} = T/\tau,
\end{equation}
where $\|\cc\| = \sqrt{\cc^T\cc}$ is the Euclidean vector norm.

In all tests we set the final time $T$ equal to $T=0.2$.
At this $T$ the solution has already passed the initial phase of 
forming homogeneous regions for times $t\approx\epsilon^2$, 
all the schemes need sufficiently many time steps, 
but the solution is still far from its stationary state.

Table~\ref{t:max_dt} presents the values of the upper spectral bound 
$\lambda_{\infty}$ depending on the grid step~$h$ and the corresponding 
maximum time step size values~$\tau$ at which the explicit scheme is stable.
In each case, two maximum $\tau$~values are given in the table:
\\ 
(a)~a maximum time step size~$\tau$, for which gradient stability
is observed (i.e., no energy increase~\eqref{psi_hn} takes place); 
(b)~a maximum time step size~$\tau$, for which the schemes
LIM-LIE and LIM-LSS operate in the explicit scheme mode,
i.e., for which relation~\eqref{p} yields $p=1$, meaning that
no Chebyshev iterations are needed.
\\
The last condition implies that the EE~scheme is stable in 
the Euclidean operator norm.  This can be seen from the following.
If for the LIM-LIE scheme relation~\eqref{p} gives $p=1$ then
\begin{equation}
\label{p=1}  
\tau\lambda_{\infty} \leqslant \left(\frac{4}{\pi}\right)^2 - 1 ,
\end{equation}
where $\lambda_{\infty} = \|\widetilde{\AB}\|_1$ with $\widetilde{\AB}$
being defined in~\eqref{lin_impl}.  
If $J_n=\partial F'(\cc^n)/\partial\cc$ then, taking into 
account approximation 
$F'(\cc^n)\approx F'(\textbf{0}) + J_n(\cc^n-\textbf{0})=J_n\cc_n$,
we obtain, for the EE~solution~$\cc^{n+1}$, 
$$
\|\cc^{n+1}\| =\|\cc^n -\tau \AB (F'(\cc^n) +\epsilon^2 \AB\cc^n)\|
\approx \|(I-\tau\widetilde{\AB})\cc^n\|\leqslant
\|I-\tau\widetilde{\AB}\| \, \|\cc^n\|. 
$$ 
Since $\widetilde{\AB}$ is symmetric positive semidefinite matrix, 
condition $\|I-\tau\widetilde{\AB}\|_2 \leqslant 1$ in the Euclidean operator norm
is equivalent to
$\tau \|\widetilde{\AB}\|_2 \leqslant 2$, 
which follows from~\eqref{p=1} 
(note that 
$\|\widetilde{\AB}\|_2\leqslant \|\widetilde{\AB}\|_1=\lambda_{\infty}$).

\begin{table}
  \centering
  \caption{The upper spectral bound $\lambda_{\infty}=\|\widehat{\AB}\|_1$
    for the LIM-LSS scheme and the maximum time step size $\tau_{\max}$, 
    for which the EE scheme is gradient stable (row~3 in the table) 
    and the LIM-LSS scheme operates in the explicit scheme mode    
    (rows~4 and~6), depending on~$h$.  
    As seen in the table, 
    $\lambda_{\infty}=\mathcal{O}(h^{-2})$ for $\epsilon=\epsilon_4(h)$,
    $\lambda_{\infty}=\mathcal{O}(h^{-4})$ for $\epsilon=\epsilon_4(1/64)$
    and in all cases $\tau_{\max}=\mathcal{O}(\lambda_{\infty}^{-1})$.}
  \label{t:max_dt}
  
  \renewcommand{\arraystretch}{1.1} 
  \begin{center}
    \begin{tabular*}{1.0\textwidth}{@{\extracolsep{\fill} } cccccc}
      \hline\hline
      $h$ & $1/32$ & $1/64$ & $1/128$ & $1/256$ & $1/512$ 
      \\\hline
      \multicolumn{6}{c}{$\epsilon=\epsilon_4(h)$}
      \\
      $\lambda_{\infty}=\|\widehat{\AB}\|_1$ & $2.3\times10^{4}$ 
                   & $9.3\times10^{4}$ & $3.7\times10^{5}$ & $1.5\times10^{6}$ & $6.0\times10^{6}$ \\
      gradient stability & $8.8\times10^{-5}$    & $2.1\times10^{-5}$     & $5.4\times10^{-6}$     & $1.3\times10^{-6}$ & $3.3\times10^{-7}$   
      \\
      explicit scheme mode & $2.6\times10^{-5}$ & $6.6\times10^{-6}$     & $1.6\times10^{-6}$     & $4.1\times10^{-7}$ & $1.0\times10^{-7}$   
      \\\hline
      \multicolumn{6}{c}{$\epsilon=\epsilon_4(1/64)$}
      \\
      $\lambda_{\infty}=\|\widehat{\AB}\|_1$ & $1.2\times10^{4}$ & $9.3\times10^{4}$ & $1.1\times10^{6}$ & $1.6\times10^{7}$ & $2.5\times10^{8}$
      \\  
      explicit scheme mode & $5.1\times10^{-5}$ & $6.6\times10^{-6}$     & $5.6\times10^{-7}$     & $3.8\times10^{-8}$ & $2.4\times10^{-9}$    \\
      \hline
    \end{tabular*}
  \end{center}
\end{table}

Table~\ref{t:max_dt} shows dependence of the upper spectral bound 
$\lambda_{\infty}$
on $h$ for the LIM-LSS scheme, i.e., $\lambda_{\infty}=\|\widehat{\AB}\|_1$.
In the tests presented here in all cases the values of
$\|\widehat{\AB}\|_1$ and $\|\widetilde{\AB}\|_1$ appear to be so close
that the Chebyshev iteration numbers in the LIM-LSS and LIM-LIE
schemes turn out to be almost the same.
The values in Table~\ref{t:max_dt} are given for the both
ways to determine $\epsilon$, see~\eqref{eps_h}, \eqref{eps_fix}.
From the table data, it is easy to check that for the choice
$\epsilon = \epsilon_4(h)$ the matrix $\epsilon^2\AB$ does not depend on $h$.
Hence, taking into account~\eqref{LSS2}, we obtain a dependence
$\|\widehat{\AB}\|_1=\mathcal{O}(h^{-2})$.  
For $\epsilon$ defined in the second way (by formula~\eqref{eps_fix}),
we have $\|\widehat{\AB}\|_1=\mathcal{O}(h^{-4})$.

\begin{figure}
\center{\includegraphics[width = 0.8\textwidth]{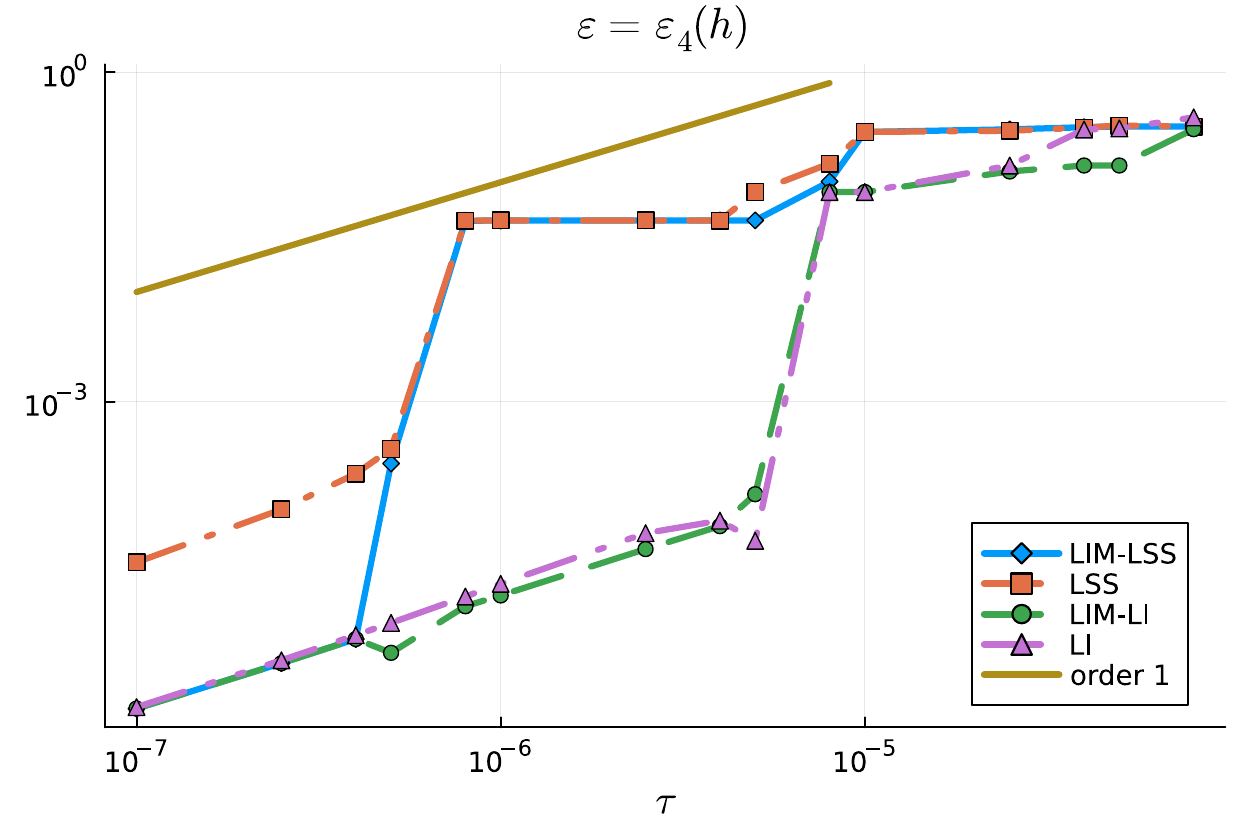}}  
\caption{Convergence (achieved accuracy versus the time step size $\tau$) 
for the LIM-LSS, LSS, LIM-LIE and LIE schemes on the $N=256$ space grid,
$\epsilon=\epsilon_4(h)$}
\label{f:conv_eps_h}
\end{figure}

The first order convergence $\mathcal{O}(\tau)$ of the schemes
LIM-LSS, LSS, LIM-LIE and LIE is confirmed by the plot 
in Figure~\ref{f:conv_eps_h}.
As can be seen from the figure, 
both LIM schemes retain the accuracy properties of the
of the schemes they are based on. 
In addition, we note that the Eyre splitting schemes LSS
and LIM-LSS turn out to be less accurate than the LIM-LIE and LIE schemes. 
In the former, the gradient stability is achieved 
by the splitting, which leads to an additional error.

In Tables~\ref{t2_ex1_N64}--\ref{t2_ex1_N512} 
computational costs and achieved accuracies,
depending on the time step size~$\tau$,
are given for
the EE, LSS, LIM-LSS, and LIM-LIE schemes.
As we see in tables, accuracies delivered by
the LSS and LIM-LSS schemes turn out to be too low
(we consider error values~\eqref{err} greater than $10^{-2}$ 
unacceptably large).  
In these schemes,
a time step size increase with respect to the EE scheme,
needed to compensate higher costs per time step,
leads to an unacceptably low accuracy.
The LIM-LIE scheme is more accurate and yields a 
gain with respect to the EE scheme on fine grids 
(the gain is approximately a factor~2 on the $N=512$ grid).
The fact that the low accuracy of the LSS and LIM-LSS is caused
by the Eyre splitting, is confirmed by switching to the
LIE and LIM-LIE schemes (the LIE scheme error values
are not shown in the table, they are close to that
of the LIM-LIE scheme, see Figure~\ref{f:conv_eps_h}).
 
\begin{table}
\caption{Number of matrix-vector multiplications (matvecs),
number of linear system solutions (for implicit schemes LSS and LIE),
and achieved accuracy versus~$\tau$.
$\epsilon = \epsilon_4(h)$, space grid $N=64$}
\label{t2_ex1_N64}
\begin{center}
    \begin{tabular*}{1.0\textwidth}{@{\extracolsep{\fill} } cccc}
        \hline\hline
        $\tau$ & scheme  & \#{} matvecs $/$ & error           \\
               &        & \#{} lin.syst.~solutions $\times10^3$  & \eqref{err} \\
        \hline
            $5.0\times10^{-5}$
            & LIM-LSS     & 12 / ---   & $5.24\times10^{-2}$     \\                
            & LSS         & 4 / 4      & $1.18\times10^{-1}$     \\
            & LIM-LIE      & 12 / ---   & $8.34\times10^{-3}$     \\ \hline
            $1.0\times10^{-5}$ 
            & EE      & 20 / ---   & $1.04\times10^{-4}$     \\
            & LIM-LSS    & 60 / ---   & $2.52\times10^{-2}$     \\
            & LSS        & 20 / 20    & $2.92\times10^{-2}$     \\
            & LIM-LIE     & 60 / ---   & $3.96\times10^{-4}$     \\\hline
            $1.0\times10^{-6}$  
            & EE      & 200 / ---~ & $2.27\times10^{-5}$     \\
            & LIM-LSS    & 200 / ---~ & $2.27\times10^{-5}$     \\
            & LSS        & 200 / 200  & $2.93\times10^{-3}$     \\
            & LIM-LIE     & 200 / ---  & $2.27\times10^{-5}$     \\
        \hline
    \end{tabular*}
\end{center}
\end{table}
 
\begin{table}
\caption{Number of matrix-vector multiplications (matvecs),
number of linear system solutions (for implicit schemes LSS and LIE),
and achieved accuracy versus~$\tau$.
$\epsilon = \epsilon_4(h)$, space grid $N=128$.}
\label{t2_ex1_N128}
\begin{center}
    \begin{tabular*}{1.0\textwidth}{@{\extracolsep{\fill} } cccc}
        \hline\hline
        $\tau$ & scheme  & \#{} matvecs $/$ & error           \\
                   &        & \#{} lin.syst.~solutions $\times10^3$  & \eqref{err} \\
        \hline
            $2.5\times10^{-5}$
            & LIM-LSS & 40 / ---   & $1.71\times10^{-1}$     \\                
            & LSS     & 8 / 8      & $1.66\times10^{-1}$     \\
            & LIM-LIE  & 40 / ---  & $9.06\times10^{-2}$     \\ \hline
            $5.0\times10^{-6}$ 
            & EE   & 40 / ---   & $3.37\times10^{-4}$     \\
            & LIM-LSS & 120 / ---~ & $3.60\times10^{-2}$ \\
            & LSS     & 40 / 40    & $1.19\times10^{-2}$     \\
            & LIM-LIE  & 120 / ---  & $7.91\times10^{-4}$     \\
        \hline
    \end{tabular*}
\end{center}
\end{table}

\begin{table}
\caption{Number of matrix-vector multiplications (matvecs),
number of linear system solutions (for implicit schemes LSS and LIE),
and achieved accuracy versus~$\tau$.
$\epsilon = \epsilon_4(h)$, space grid $N=256$.}
\label{t2_ex1_N256}
\begin{center}
    \begin{tabular*}{1.0\textwidth}{@{\extracolsep{\fill} } cccc}
        \hline\hline
        $\tau$ & scheme  & \#{} matvecs $/$ & error           \\
               &        & \#{} lin.syst.~solutions $\times10^3$  & \eqref{err} \\
        \hline
$1.0\times10^{-5}$ & LIM-LIE  & 140 / ---  & $8.14\times10^{-2}$     \\ \hline
$5.0\times10^{-6}$ & LIM-LSS & 200 / ---  & $4.49\times10^{-2}$     \\             
                  & LSS     & 40 / 40    & $8.12\times10^{-2}$     \\ 
                  & LIM-LIE  & 200 / ---  & $1.44\times10^{-4}$     \\ \hline
$1.0\times10^{-6}$ 
            & EE   & 200 / ---    & $2.06\times10^{-5}$ \\
            & LIM-LSS & 600 / ---    & $4.48\times10^{-2}$ \\
            & LSS     & 200 / 200    & $4.48\times10^{-2}$     \\
            & LIM-LIE  & 600 / ---    & $1.72\times10^{-5}$     \\ 
        \hline
    \end{tabular*}
\end{center}
\end{table}

\begin{table}
\caption{Number of matrix-vector multiplications (matvecs),
number of linear system solutions (for implicit schemes LSS and LIE),
and achieved accuracy versus~$\tau$,
$\epsilon = \epsilon_4(h)$, space grid $N=512$.}
    \label{t2_ex1_N512}
	\begin{center}
		\begin{tabular*}{1.0\textwidth}{@{\extracolsep{\fill} } cccc}
			\hline\hline
			$\tau$ & scheme  & \#{} matvecs $/$ & error           \\
                       &        & \#{} lin.syst.~solutions $\times10^6$  & \eqref{err} \\
			\hline
$2.0\times10^{-6}$ 
                & LIM-LIE  & 0.5 / ---  & $8.25\times10^{-6}$     \\\hline
                $1.0\times10^{-6}$
                & LIM-LSS & 1 / ---   & $3.25\times10^{-2}$     \\                
                & LSS     & 0.2 / 0.2      & $3.25\times10^{-2}$     \\
                & LIM-LIE  & 1 / ---  & $4.61\times10^{-6}$     \\ \hline
$2.0\times10^{-7}$ 
                & EE   & 1 / ---   & $9.57\times10^{-7}$    \\
                & LIM-LSS & 3 / ---~ & $3.25\times10^{-2}$     \\
                & LSS     & 1 / 1    & $3.25\times10^{-2}$     \\
                & LIM-LIE  & 3 / ---  & $7.35\times10^{-7}$     \\ 
			\hline
		\end{tabular*}
	\end{center}
 \end{table}

\begin{figure}
\center{\includegraphics[width = 0.8\textwidth]{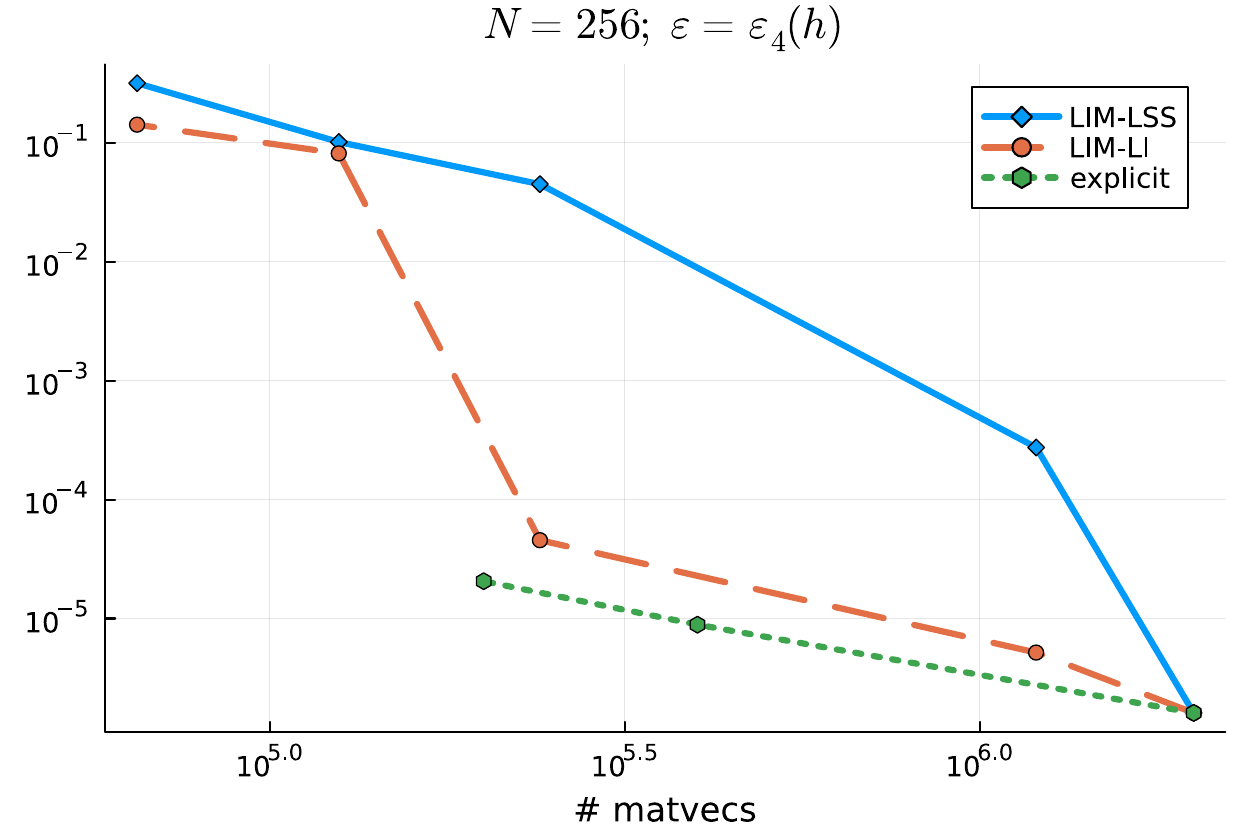}}  
\center{\includegraphics[width = 0.8\textwidth]{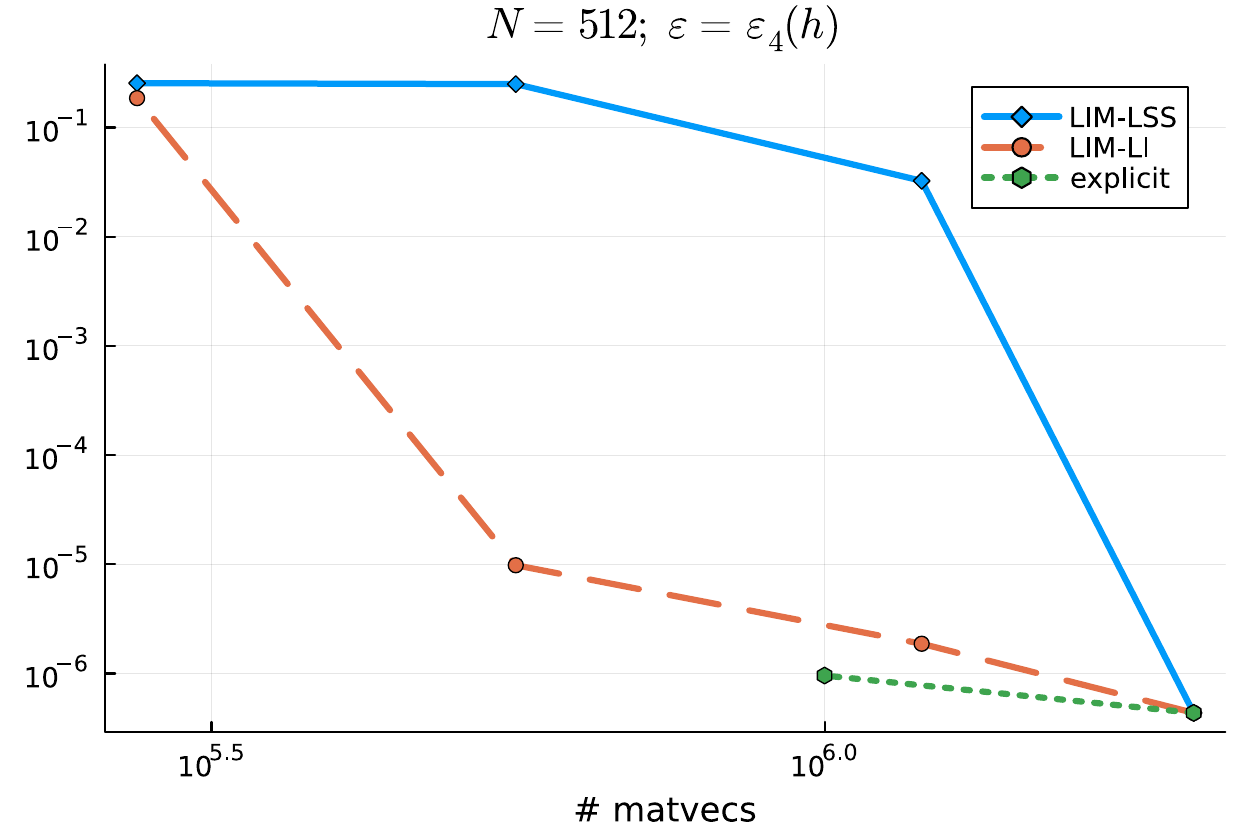}}  
\caption{Achieved accuracy versus number of matrix-vector multiplications
(matvecs) for the EE and LIM schemes on the 
$N=256$ space grid (upper plot) and
$N=512$ space grid (lower plot), $\epsilon=\epsilon_4(h)$.
Increasing $\tau$ in the EE scheme further is impossible
due to the stability restrictions.}
\label{f:err_wrk1}  
\end{figure}

Plots showing achieved accuracy versus
computational costs are given in Figure~\ref{f:err_wrk1}.
As can be seen in the plots, the LIM schemes perform
better on the $N=512$ grid than on the $N=256$ grid:
a cost reduction of an approximately factor~3 is achieved
on the $N=512$ grid.

In the presented tests a rather low achieved accuracy is observed
for all the schemes except the EE~scheme and this is not only due 
to the Eyre splitting.  In our tests the initial value vector
is a non-smooth grid function, where each vector entry is chosen
randomly.
As discussed above in the introduction, a fast forming 
of the homogeneity regions is observed for such initial values,
taking place at typical evolution times $\epsilon^2$.  
Hence, to properly track this process, the time step size
should be chosen as $\tau \sim \mathcal{O}(\epsilon^2)= \mathcal{O}(h^2)$.
Therefore, taking into account that 
$\epsilon^2\|A\|\sim \epsilon^2h^{-2} = \mathcal{O}(1)$,
we see that during this initial phase of the time integration
the time step size~$\tau$ can not be taken much larger than
the time step size in the EE scheme.
This time step is, in this case, determined by the accuracy rather than
by the stability requirements.

To test the stability and accuracy of the considered time integration schemes
at larger time steps, we carry out tests with a smoother initial condition. 
For the space grids~$N>64$, we set it by interpolating the initial non-smooth 
initial solution from the grid $N = 64$ to the current (finer) space grid. 
On the $N = 64$ space grid the initial vector does not change. 
The interpolation is done by the piecewise cubic Hermite interpolation, 
which excludes the occurrence of new extrema and yields 
a continuously differentiable function (in the octave package, 
this interpolation method is called \texttt{pchip}). 
In addition, in the tests with the smoothed initial value vector, 
the second method~\eqref{eps_fix} for setting $\epsilon$ is employed, 
which allows to test our schemes in the situation
when the norm of the right-hand side operator grows as $\mathcal{O}(h^{-4})$,
see Table~\ref{t:max_dt}.

\begin{figure}
\center{\includegraphics[width = 0.8\textwidth]{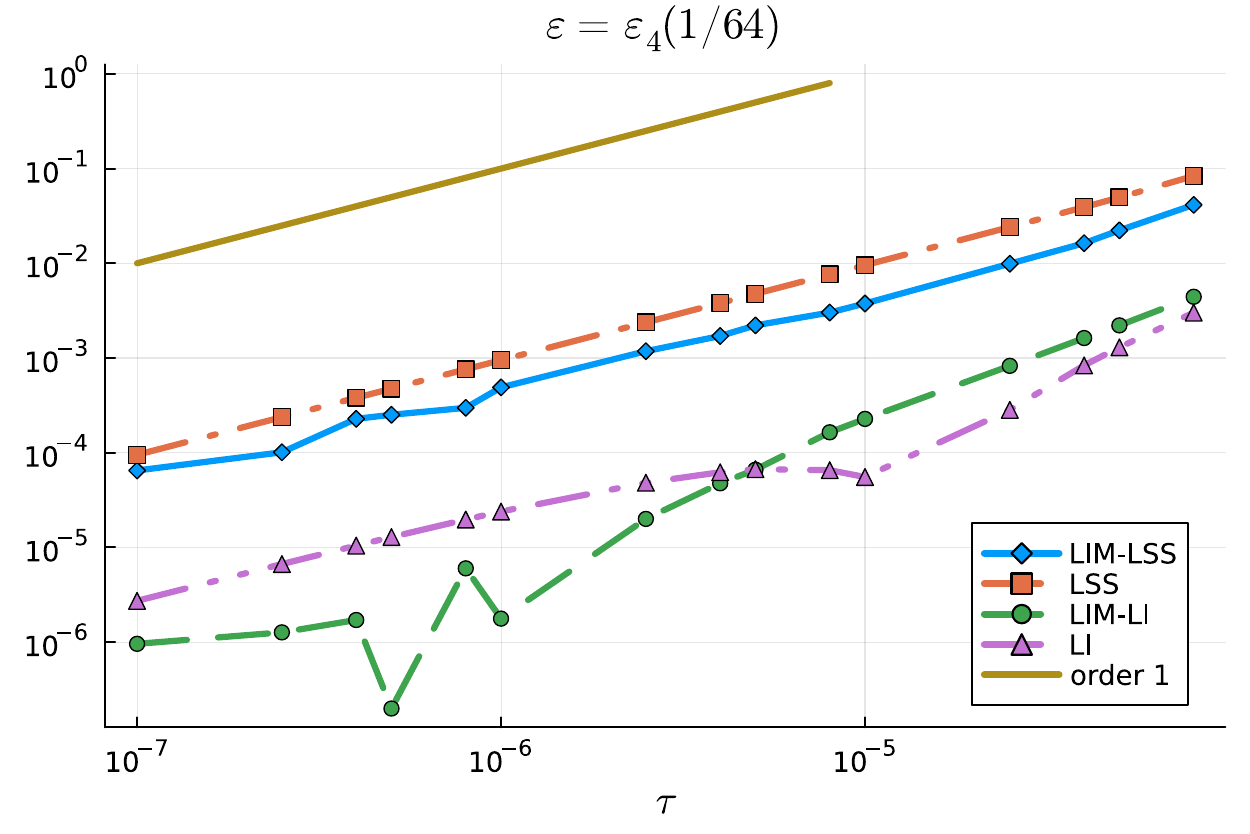}} 
\caption{Convergence (achieved accuracy versus the time step size $\tau$) 
for the LIM-LSS, LSS, LIM-LIE and LIE schemes on the $N=256$ space grid,
$\epsilon=\epsilon_4(1/64)$}
\label{f:conv_eps_fix}
\end{figure}

We start by testing the first order convergence $\mathcal{O}(\tau)$ 
at the plot in Figure~\ref{f:conv_eps_fix}.
Comparing the plot with the plot in~\ref{f:conv_eps_h}, 
we see that all the schemes achieve a much higher accuracy.
However, the LSS and LIM-LSS schemes based on the Eyre splitting 
are still less accurate than the LIE and LIM-LIE schemes.
Small error oscillations observed with the LIM-LIE scheme 
for $\tau\approx 10^{-6}$ are caused by the fact that 
number of Chebyshev iterations carried out each time step
slightly varies for these~$\tau$ values. 

\begin{figure}
\center{\includegraphics[width = 0.8\textwidth]{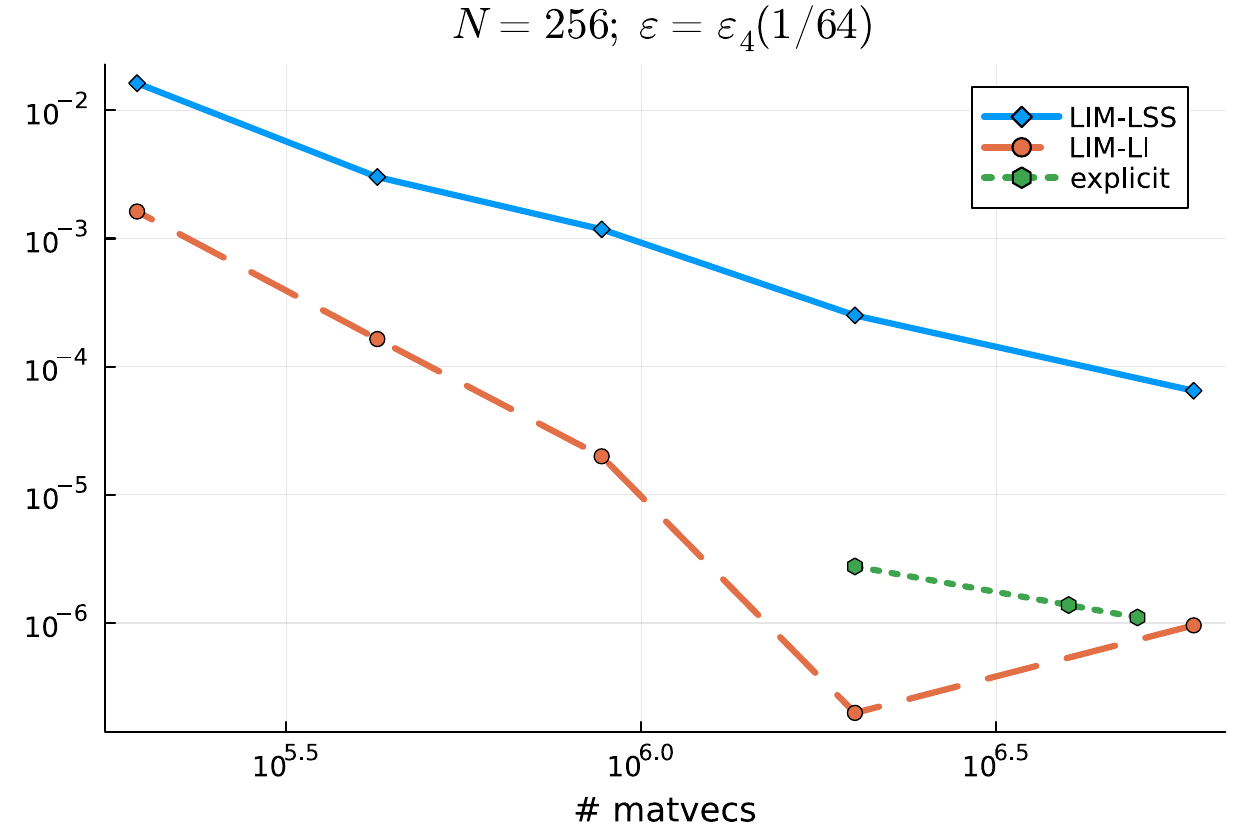}}  
\center{\includegraphics[width = 0.8\textwidth]{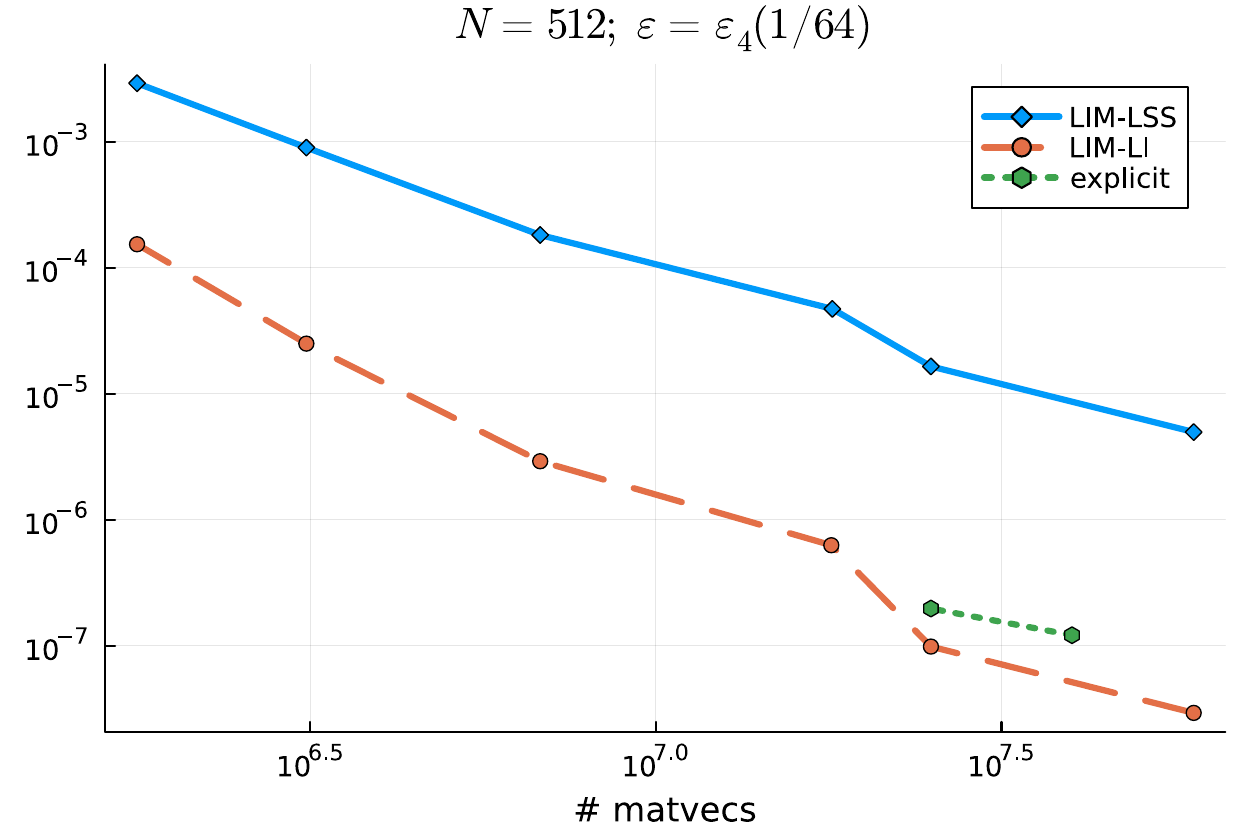}}  
\caption{Achieved accuracy versus the number of matrix-vector multiplications
(matvecs) for the EE scheme and the LIM schemes
with the smoothed initial vector~$\cc^0$ and $\epsilon=\epsilon_4(1/64)$
on the $N=256$ grid (upper plot) and the $N=512$ grid (bottom plot).
Increasing $\tau$ in the EE scheme further is impossible
due to the stability restrictions.}
\label{f:err_wrk2}  
\end{figure}

Test results for smoothed initial value vectors and fixed
$\epsilon=\epsilon_4(1/64)$ are presented in 
Tables~\ref{t2_ex1_N128_sm}--\ref{t2_ex1_N512_sm} and 
in Figure~\ref{f:err_wrk2}.  We see that both LIM schemes
provide an essential efficiency gain with respect to the EE~scheme.  
On the $N=512$ grid the gain of approximately a factor~10
is attained (see the bottom plot in Figure~\ref{f:err_wrk2}).

\begin{table}
\caption{Number of matrix-vector multiplications (matvecs),
number of linear system solutions (for implicit schemes LSS and LIE),
and achieved accuracy versus~$\tau$.
$\epsilon = \epsilon_4(h)$, smoothed initial vector~$\cc^0$, 
space grid~$N=128$.}
\label{t2_ex1_N128_sm}
	\begin{center}
		\begin{tabular*}{1.0\textwidth}{@{\extracolsep{\fill} } cccc}
			\hline\hline
			$\tau$ & scheme  & \#{} matvecs $/$ & error           \\
                       &        & \#{} lin.syst.~solutions $\times10^3$  & \eqref{err} \\
			\hline
$1.0\times10^{-4}$ 
                & LIM-LSS & 34 / ---    & $7.46\times10^{-2}$     \\
                & LSS     &  2 / 2      & $1.37\times10^{-1}$     \\ 
                & LIM-LIE  & 34 / ---    & $6.06\times10^{-3}$     \\ \hline
$5.0\times10^{-6}$
                & LIM-LSS & 200 / ---   & $3.42\times10^{-3}$     \\ 
                & LSS     & 40 / 40     & $5.25\times10^{-3}$     \\
                & LIM-LIE & 185 / ---    & $1.09\times10^{-4}$     \\ \hline
$1.0\times10^{-6}$ 
                & EE   & 200 / ---  & $3.02\times10^{-5}$     \\
                & LIM-LSS & 600 / ---~ & $8.14\times10^{-4}$     \\
                & LSS     & 200 / 200  & $1.04\times10^{-3}$     \\
                & LIM-LIE  & 600 / ---  & $1.50\times10^{-5}$     \\ 
			\hline
		\end{tabular*}
	\end{center}
 \end{table}
 
\begin{table}
\caption{Number of matrix-vector multiplications (matvecs),
number of linear system solutions (for implicit schemes LSS and LIE),
and achieved accuracy versus~$\tau$.
$\epsilon = \epsilon_4(h)$, smoothed initial vector~$\cc^0$, 
space grid~$N=256$.}
\label{t2_ex1_N256_sm}
\begin{center}
		\begin{tabular*}{1.0\textwidth}{@{\extracolsep{\fill} } cccc}
			\hline\hline
			$\tau$ & scheme  & \#{} matvecs $/$ & error           \\
                       &        & \#{} lin.syst.~solutions $\times10^6$  & \eqref{err} \\
			\hline
$1.0\times10^{-5}$ 
                & LIM-LSS & 0.38 / ---  & $3.76\times10^{-3}$     \\ 
                & LSS     & 0.02 / 0.02      & $9.55\times10^{-3}$     \\ 
                & LIM-LIE  & 0.38 / ---  & $2.27\times10^{-4}$     \\ \hline
$5.0\times10^{-7}$
                & LIM-LSS & 2 / ---   & $2.51\times10^{-4}$     \\ 
                & LSS     & 0.4 / 0.4      & $4.75\times10^{-4}$     \\
                & LIM-LIE  & 2 / ---   & $1.99\times10^{-7}$     \\ \hline
$1.0\times10^{-7}$ 
                & EE   & 2 / ---   & $2.76\times10^{-6}$     \\
                & LIM-LSS & 6 / ---   & $6.49\times10^{-5}$     \\
                & LSS     & 2 / 2     & $9.49\times10^{-5}$     \\
                & LIM-LIE  & 6 / ---   & $9.59\times10^{-7}$     \\ 
			\hline
		\end{tabular*}
	\end{center}
 \end{table}
 
\begin{table}
\caption{Number of matrix-vector multiplications (matvecs),
number of linear system solutions (for implicit schemes LSS and LIE),
and achieved accuracy versus~$\tau$.
$\epsilon = \epsilon_4(h)$, smoothed initial vector~$\cc^0$, 
space grid~$N=512$.}
\label{t2_ex1_N512_sm}
\begin{center}
\begin{tabular*}{1.0\textwidth}{@{\extracolsep{\fill} } cccc}
  \hline\hline
  $\tau$ & scheme  & \#{} matvecs $/$ & error           \\
         &        & \#{} lin.syst.~solutions $\times10^6$  & \eqref{err} \\
  \hline
  $8.0\times10^{-7}$ 
  & LIM-LSS & 5.75 / ---  & $3.14\times10^{-4}$     \\
  & LSS     & 0.25 / 0.25   & $7.13\times10^{-4}$     \\ 
  & LIM-LIE  & 5.75 / ---  & $3.83\times10^{-6}$     \\ \hline
  $4.0\times10^{-8}$
  & LIM-LSS & 25 / ---   & $1.65\times10^{-5}$     \\ 
  & LSS     &  5 / 5      & $3.56\times10^{-5}$     \\
  & LIM-LIE  & 25 / ---  & $9.87\times10^{-8}$     \\ \hline
  $8.0\times10^{-9}$ 
  & EE   & 25 / ---   & $1.98\times10^{-7}$    \\
  & LIM-LSS & 75 / ---~ & $4.44\times10^{-6}$     \\
  & LSS     & 25 / 25    & $7.13\times10^{-6}$     \\
  & LIM-LIE  & 75 / ---  & $5.14\times10^{-8}$     \\ 
  \hline
\end{tabular*}
\end{center}
\end{table}

\section{Conclusions}\label{sec:conclusions}

The presented results allow to make the following
conclusions.

\begin{enumerate}
\item 
Proposed local iteration schemes LIM-LSS and LIM-LIE~\eqref{LI-M},
which are based on the implicit schemes
LSS~\eqref{LSS} and LIE~\eqref{lin_impl},
proved to be reliable in practice.
They combine a simplicity and parallelism of explicit schemes
with stability of implicit schemes.
Theoretical estimates of efficiency for the LIM schemes
are confirmed in the tests: the LIM schemes provide 
an efficiency gain up to a factor~10 with respect to 
the explicit scheme~EE.
The gain factor grows as the space grid gets finer.

\item
In the considered numerical tests,
gradient stable schemes based on the Eyre splitting appear
to be less accurate than regular linearized schemes.
In particular, the regular linearized implicit Euler scheme~LIE~\eqref{lin_impl}
and based on it local iteration scheme LIM-LIE provide a higher
accuracy than the linear stabilized splitting scheme~LSS~\eqref{LSS} 
and the local iteration scheme LIM-LSS.
In the latter schemes a gradient stability is achieved by splitting,
which leads to an additional error.

\item
In cases where the diffusion boundary width $\epsilon$ is chosen proportional
to the grid size, i.e., $\epsilon=\mathcal{O} (h)$, 
the time step size in the explicit scheme~EE is bounded as $\tau=\mathcal{O}(h^2)$.
Therefore, we can expect that the potential of the local
iteration schemes for the Cahn-Hilliard equation is comparable
to that for parabolic problems.

\item
For non-smooth initial data, for instance, if the initial value vector
is chosen randomly, additional accuracy restrictions on the time step
size~$\tau$ arise at times $t \leqslant\mathcal{O}(\epsilon^2)$.  
Since decreasing the time step size in an implicit scheme 
usually means lowering its efficiency,
the local iteration schemes appear to be especially attractive
(in these schemes decreasing the time step size does lead to lower 
computational costs).
In such problems it seems sensible to apply local iterations
schemes with an adaptive choice of the time step size~\cite{BotZhukov2024}.   
\end{enumerate}

\subsection*{Acknowledgments}
The authors thank V.T.~Zhukov (KIAM RAS) for useful discussions
and consulting on local iterations schemes.

\end{document}